\documentclass[12pt]{amsart}
\usepackage{a4wide,amsfonts,amssymb,stmaryrd,amscd,amsmath,latexsym,amsbsy}
\numberwithin{equation}{section}
\theoremstyle{plain}
\newtheorem{thm}{Theorem}[section]

\newtheorem{prop}[thm]{Proposition}

\newtheorem{defi}[thm]{Definition}
\newtheorem{lem}[thm]{Lemma}
\newtheorem{cor}[thm]{Corollary}
\newtheorem{eg}[thm]{{Example}}

\theoremstyle{remark}
\newtheorem{rema}[thm]{Remark}

\newcommand {\CQG}{\textsf{CQG~}} 

\newcommand{\Fb}{\mathbb{F}}      
\newcommand {\Fz}{\mathbb{E}}     

\newcommand{\Ko}{\mathcal{K}}      
\newcommand{\Le}{\mathcal{M}}      
\newcommand{\KF}{\Ko_{\Fb}}       
\newcommand{\LF}{\Le_{\Fb}}       

\newcommand{\Gr}{\mathcal{G}}     
\newcommand{\RR}{\mathbb{C}}     

\newcommand{\xx}{x}     
\newcommand{\yy}{y}     
\newcommand{\kk}{K}     
\newcommand{\zz}{X}     

\title{Quantum dimensions and their non-Archimedean degenerations}
\author{Uri Onn}
\address{Uri Onn, Einstein Institute of Mathematics, Edmond Safra Campus, Givat Ram,
Jerusalem 91904, Israel.}
\email{urion@math.huji.ac.il}
\author{Jasper V. Stokman}
\address{J.V. Stokman,
KdV Institute for Mathematics, Universiteit van Amsterdam,
Plantage Muidergracht 24, 1018 TV Amsterdam, The Netherlands.}
\email{jstokman@science.uva.nl}

\subjclass[2000]{33D52, 33D80}

\thanks{The first author is supported by the Israel Science
Foundation, ISF grant no. 555104, and by the Edmund Landau Minerva
Center for Research in Mathematical Analysis and Related Areas,
sponsored by the Minerva Foundation (Germany). The second author
is supported by the Netherlands Organization for Scientific Research
(NWO) in the VIDI-project  \lq Symmetry and modularity in exactly
solvable models\rq.}

\begin{document}

\begin{abstract}
We derive explicit dimension formulas for irreducible
$\mathcal{M}_{\mathbb{F}}$-spherical
$\mathcal{K}_{\mathbb{F}}$-re\-pre\-sen\-ta\-tions where
$\mathcal{K}_{\mathbb{F}}$ is the maximal compact subgroup of the
general linear group $\hbox{GL}_d(\mathbb{F})$ over a local field
$\mathbb{F}$ and $\mathcal{M}_{\mathbb{F}}$ is a closed subgroup of
$\mathcal{K}_{\mathbb{F}}$ such that
$\mathcal{K}_{\mathbb{F}}/\mathcal{M}_{\mathbb{F}}$ realizes the
Grassmannian of $n$-dimensional $\mathbb{F}$-subspaces of
$\mathbb{F}^d$. We explore the fact that $(\mathcal{K}_{\mathbb{F}},
\mathcal{M}_{\mathbb{F}})$ is a Gelfand pair whose associated zonal
spherical functions identify with various degenerations of the
multivariable little $q$-Jacobi polynomials. As a result, we are led
to consider generalized dimensions defined in terms of evaluations
and quadratic norms of multivariable little $q$-Jacobi polynomials,
which interpolate between the various classical dimensions. The
generalized dimensions themselves are shown to have representation
theoretic interpretations as the quantum dimensions of irreducible
spherical quantum representations associated to quantum complex
Grassmannians.
\end{abstract}

\maketitle



\section{Introduction}
Let $\Fb$ be a local field and
let $\KF$ be the maximal compact subgroup of $\hbox{GL}_d(\Fb)$.
That is,
\begin{equation}\label{KF}
\KF=
\begin{cases}
\hbox{O}(d)\qquad &\hbox{ if } \Fb=\mathbb{R},\\
\hbox{U}(d)\qquad &\hbox{ if } \Fb=\mathbb{C},\\
\hbox{GL}_d(\mathcal{O})\qquad &\hbox{ if } \Fb \,\,
\hbox{non-Archimedean,}
\end{cases}
\end{equation}
where in the latter case $\mathcal{O}$ denotes the ring of integers
of $\Fb$. Let $\Gr(n,d;\Fb)$ be the Grassmannian of $n$-dimensional
subspaces of a fixed $d$-dimensional vector space over $\Fb$.
Throughout we assume that $n\leq d/2 $. The group $\KF$ acts
transitively on $\Gr(n,d;\Fb)$. Letting $\LF$ be a stabilizer of a
point we may identify $\Gr(n,d;\Fb)$ with $\KF/\LF$. The pair
$(\KF,\LF)$ is a Gelfand pair, in the sense that each (continuous,
complex) irreducible $\KF$-representation has a subspace of
$\LF$-fixed vectors which is at most one-dimensional. The
irreducible representations having a one-dimensional subspace of
$\LF$-fixed vectors are called $\LF$-spherical. In this paper we derive
explicit dimension formulas for the irreducible $\LF$-spherical
$\KF$-representations.

For all local fields $\Fb$ the equivalence classes of the
irreducible $\LF$-spherical $\KF$-re\-pre\-sen\-ta\-tions are
naturally parameterized by the set $\Lambda_n$ of partitions of at
most $n$ parts (see e.g. \cite{JC} and \cite{DS} for $\Fb$
Archimedean, and \cite{BO2} for $\Fb$ non-Archimedean). This is
quite remarkable, taking into consideration that the representation
theory of $\KF$ for $\Fb$ non-Archimedean is of a completely
different nature compared to the representation theory for $\Fb$
Archimedean. In each case there is a natural explicit labeling of
the irreducible $\LF$-spherical $\KF$-representations
$V_\lambda^{\Fb}$ by partitions $\lambda\in\Lambda_n$ which we shall
describe in detail.

The zonal spherical functions $\varphi_\lambda^{\Fb}$
associated to $V_\lambda^{\Fb}$ ($\lambda\in\Lambda_n$)
form a linear basis of the Hecke algebra of $\LF$-biinvariant
representative functions on $\KF$. The starting point for the derivation of the
explicit dimension formulas for the irreducible $\LF$-spherical $\KF$-representations
is the representation theoretic dimension formula
\begin{equation}\label{repdimintro}
\textup{Dim}_{\mathbb{C}}\bigl(V_\lambda^{\Fb}\bigr)=
\frac{\bigl(\varphi_\lambda^{\Fb}(e)\bigr)^2}
{\int_{\KF}|\varphi_\lambda^{\Fb}(g)|^2dg},\qquad
\lambda\in\Lambda_n
\end{equation}
where
$e\in\KF$ is the unit element and integration is with respect to
the normalized Haar measure on $\KF$. The second key ingredient is
the fact that the zonal spherical functions $\varphi_\lambda^{\Fb}$
for the different local fields $\Fb$ relate to various
degenerations of the multivariable little $q$-Jacobi polynomials
from \cite{S0}, see \cite{JC}, \cite{DS} and \cite{OU1}.
As a result, we are led to consider generalized dimensions
defined in terms of evaluations and quadratic norms of multivariable little $q$-Jacobi
polynomials, which interpolate between the various classical
dimensions.

The multivariable little $q$-Jacobi polynomials $P_\lambda^L(z)$ ($\lambda\in\Lambda_n$)
are symmetric polynomials in $n$ variables $z=(z_1,\ldots,z_n)$
depending on four auxiliary parameters, which arise
as degenerations of the Macdonald-Koornwinder polynomials \cite{K}. They form an orthogonal system with respect
to an explicit discrete probability orthogonality measure
defined in terms of an iterated Jackson integral (see \cite{S0}),
which induces a norm $\|\cdot\|_L$ on the symmetric polynomials in
$z$. The generalized dimensions referred to in the previous paragraph are
\begin{equation}\label{tildeD}
d(\lambda)=\frac{P_\lambda^L(0)^2}{\|P_\lambda^L\|_L^2},\qquad
\lambda\in\Lambda_n,
\end{equation}
compare with the representation theoretic dimension formula \eqref{repdimintro}.

We show that the dimension of $V_\lambda^{\Fb}$
can be obtained as a classical ($q=1$) degeneration of
$d(\lambda)$ for $\Fb$ Archimedean, and as a $p$-adic
($q=0$) degeneration of $d(\lambda)$ for $\Fb$
non-Archimedean. For special values of the four auxiliary parameters
of the multivariable little $q$-Jacobi polynomials we
show that the generalized dimensions $d(\lambda)$
($\lambda\in\Lambda_n$) themselves have a representation theoretic interpretation
as the quantum dimensions of irreducible spherical quantum representations
associated to the one-parameter family of quantum complex
Grassmannians from \cite{NDS} and \cite{DS}.
In this case the key ingredient is the identification in \cite{NDS}
and \cite{DS} of the associated quantum zonal spherical functions
with Macdonald-Koornwinder polynomials, as well as with degenerations
of the Macdonald-Koornwinder polynomials known as
multivariable big and little $q$-Jacobi polynomials.

Cherednik's double affine Hecke algebra techniques have led to
explicit evaluation formulas and to explicit quadratic norm evaluations for the
Macdonald-Koornwinder polynomials in \cite{Sa} (see also \cite{vD}
and \cite{S2}). The quadratic norms $\|P_\lambda\|_L^2$ of the multivariable
little $q$-Jacobi polynomials (as well as for multivariable big $q$-Jacobi polynomials)
have been explicitly evaluated in \cite{S} by degenerating the quadratic norm evaluations of
the Macdonald-Koornwinder polynomials. In exactly the same way we
derive in this paper evaluation formulas for the multivariable big
and little $q$-Jacobi polynomials. This leads to
explicit expressions for the interpolating, generalized dimensions $d(\lambda)$
($\lambda\in\Lambda_n$). Explicit (quantum) dimension formulas for the irreducible
spherical (quantum) representations associated to the various (quantum)
Grassmannians are subsequently derived by degenerating the
explicit expression for $d(\lambda)$ ($\lambda\in\Lambda_n$).

The contents of the paper is as follows. In Section
\ref{qJsection} we introduce the multivariable little $q$-Jacobi
polynomials and the associated generalized dimensions. In Section
\ref{CQG} we discuss the representation theory of general Compact
Quantum Group (\textsf{CQG}) algebras from \cite{DK} and \cite{F}.
\CQG algebras form the natural setting to capture the harmonic
analytic structures of the various (quantum) Grassmannians under
consideration in this paper.

In Section \ref{local} we specialize the general representation
theoretic setup from Section \ref{CQG} to the classical \CQG algebra
of representative functions on the compact group $\KF$.
We discuss the harmonic analytic implications for the associated
Grassmannian $\KF/\LF$. We pay close attention to the specific
parametrization of the irreducible $\LF$-spherical
$\KF$-representations by $\Lambda_n$. We state the explicit dimension formulas,
which are new in case of non-Archimedean local
fields $\Fb$ (see Theorem \ref{thm-dim-padic}). Although the results are
essentially classical for Archimedean local fields $\Fb$ (see e.g.
\cite{JC}, \cite{DS}, \cite{Su} and \cite{OU1}), we have decided to give a
detailed exposition in case of $\Fb=\mathbb{C}$ in order to highlight the
similarities and differences to the results for the Grassmannians over non-Archimedean
local fields $\Fb$, as well as to the results for the quantum complex
Grassmannians in the subsequent sections.

In Section \ref{qGrassmannian} we specialize the setup from Section
\ref{CQG} to the \CQG algebra of quantized representative functions
on $\Ko_{\mathbb{C}}=\hbox{U}(n)$, and we discuss the associated
harmonic analysis on the standard quantum complex Grassmannian. In
particular, we state explicit quantum dimension formulas for the
associated irreducible spherical quantum representations.

In Section \ref{evaluation} we introduce the Macdonald-Koornwinder
polynomials and the multivariable big and little $q$-Jacobi
polynomials. We derive explicit evaluation formulas for the multivariable big and
little $q$-Jacobi polynomials. Subsequently we derive in
Section \ref{Dimensionsection} explicit expressions for
the generalized dimensions as well as for their $p$-adic ($q=0$)
degenerations.

In Section \ref{quantum} we discuss the present approach for the
one-parameter family of quantum complex Grassmannians from
\cite{NDS} and \cite{DS}. Following closely the analysis
of Noumi \cite{N} for other examples of quantum symmetric spaces, we express the
quantum dimensions of the associated irreducible spherical quantum
representations in terms of evaluations and quadratic norms of
Macdonald-Koornwinder polynomials, as well as of multivariable big and little
$q$-Jacobi polynomials. In particular, we obtain the
representation theoretic interpretation of the generalized dimensions for special
parameters as quantum dimensions of irreducible spherical quantum representations.
In this context we remark that
a similar approach can be followed for quantum real Grassmannians using the
harmonic analytic results from \cite{Su} and \cite{L}, but we do
not pursue it in this paper.

\noindent
{\bf{Notations and conventions.}} Representations of topological groups are complex and
are required to be continuous.
We set
\begin{equation}\label{Fz}
\Fz=\mathbb{C}(a,b,c,d,q,t)
\end{equation}
for the field of rational functions in six indeterminates $a,b,c,d,q$ and $t$.
We set $\Lambda_n$ for the partitions of at most $n$ parts. It is the
cone
\[\Lambda_n=\bigoplus_{r=1}^n\mathbb{Z}_+\omega_r
\]
in $\mathbb{Z}^n$, where $\omega_r=(1^r,0^{n-r})\in\Lambda_n$ is the
fundamental partition consisting of $r$ ones and $n-r$ zeros. We
write $\omega_0=(0^n)\in\Lambda_n$ for the zero partition. For a
partition $\lambda\in\Lambda_n$ we denote $\lambda^\prime$ for the
conjugate partition and we set $\lambda_0^\prime=n$. In addition we
write $|\lambda|$ for the weight of $\lambda$,
$l(\lambda)=\lambda_1^\prime$ for the length of $\lambda$, and
$\partial\lambda^\prime=(\lambda_j^\prime-\lambda_{j+1}^\prime)_{j\geq
0}$ for the positive integral differences of the conjugate partition
(whose sum equals $n$). For $n$-vectors $u,v,w$ we set $uv^w=(u_1
v_1^{w_1},\ldots,u_n v_n^{w_n})$. We occasionally use this notation
when $u$ or $v$ are scalars, in which we use the corresponding
$n$-vectors having all entries equal to the given scalar. We recall
some standard notations from basic hypergeometric series (see e.g.
\cite{GR}). For $j\in\mathbb{Z}_{+}$ we set
$\bigl(a;q\bigr)_j=(1-a)(1-aq)\cdots (1-aq^{j-1})$ (empty product is
one) for the $q$-shifted factorial. We write
$\bigl(a_1,\ldots,a_m;q\bigr)_j=\prod_{s=1}^m\bigl(a_s;q\bigr)_j$
for products of $q$-shifted factorials. For
$m,l,l_1,\ldots,l_k\in\mathbb{Z}_+$ with $l\leq m$, we use the
notations
\begin{equation*}
\begin{split}
\lbrack m\rbrack_{}&=\frac{1-q^{m}}{1-q},\qquad \qquad\qquad\lbrack
m\rbrack_{q}!=\prod_{j=1}^m\lbrack
j\rbrack_{q},\\
\bigg[{m \atop l_1~ l_2~ \cdots ~l_k}\bigg]_{q}&=
\frac{[m]_{q}!}{[l_1]_{q}! \cdots [l_k]_{q}!}, \qquad \left({m \atop
l}\right)_{q}=\bigg[{m \atop l~~ m-l}\bigg]_{q}.
\end{split}
\end{equation*}
In terms of $q$-shifted factorials, the $q$-factorial can be expressed as
\begin{equation}\label{brackq}
\lbrack
m\rbrack_{q}!=q^{\frac{1}{2}m(m-1)}\frac{\bigl(q^{-1};q^{-1}\bigr)_m}{(1-q^{-1})^m}.
\end{equation}
We use standard notations for the basic hypergeometric series
${}_{r+1}\phi_r$, see e.g. \cite{GR}.



\section{Multivariable little $q$-Jacobi
polynomials}\label{qJsection}

Important in our approach is the fact that (quantum) dimensions of
irreducible spherical representations associated to quantum,
real/complex and $p$-adic Grassmannians are degenerations of a
generalized dimension formula involving multivariable little
$q$-Jacobi polynomials, which we now introduce. We work over the
field $\Fz$ (see \eqref{Fz}), although in the present situation $c$
and $d$ are dummy parameters and could just as well be omitted (in
contrast to the Macdonald-Koornwinder and multivariable big
$q$-Jacobi case as treated in Section \ref{evaluation}).

The symmetric group $S_n$ in $n$ letters acts on $\mathbb{Z}_+^n$ by
permuting the coordinates. The corresponding fundamental domain is
the cone $\Lambda_n$ of partitions with at most $n$ parts.
The symmetric group $S_n$ acts on
$\Fz[z]=\Fz[z_1,\ldots,z_n]$ by permuting the independent variables
$z_1,\ldots,z_n$. We denote $\Fz[z]^S$ for the subalgebra of
$S_n$-invariant polynomials in $\Fz[z]$. An $\Fz$-basis of
$\Fz[z]^S$ is formed by the symmetric monomials
$\widetilde{m}_\lambda(z)=\sum_{\mu\in S_n\lambda}z^\mu$
($\lambda\in\Lambda_n$). The monic multivariable little $q$-Jacobi
polynomial $P_\lambda^L(z)=P_\lambda^L(z;a,b;q,t)\in \Fz[z]^S$ of
degree $\lambda\in\Lambda_n$ from \cite{S0} is of the form
\begin{equation}\label{expansionL}
P_\lambda^L(z)=\widetilde{m}_\lambda(z)+\sum_{\mu\in\Lambda_n: \mu<\lambda}c_{\lambda,\mu}^L\widetilde{m}_{\mu}(z)
\end{equation}
for certain coefficients
$c_{\lambda,\mu}^L=c_{\lambda,\mu}^L(a,b;q,t)\in \Fz$. The multivariable little $q$-Jacobi
polynomials can be characterized as solution of a second-order
difference equation, or in terms of orthogonality relations defined
with respect to an explicit discrete probability orthogonality
measure. The characterization in terms of the difference equation is
as follows. For $j=1,\ldots,n$ we write $T_j$ for the multiplicative
$q$-shift in $z_j$,
\[(T_jf)(z)=f(z_1,\ldots,z_{j-1},qz_j, z_{j+1},\ldots,z_n),\qquad f\in \Fz[z^{\pm 1}].
\]
The characterizing difference equation then reads $D_LP_\lambda^L=E_\lambda^{qJ}P_\lambda^L$ where
\begin{equation*}
\begin{split}
D_L&=\sum_{j=1}^n\bigl(\phi_{L,j}^+(z)(T_j-\hbox{Id})+\phi_{L,j}^-(z)(T_j^{-1}-\hbox{Id})\bigr),\\
\phi_{L,j}^+(z)&=qt^{n-1}a\bigl(b-\frac{1}{qz_j}\bigr)\prod_{l\not=j}\frac{z_l-tz_j}{z_l-z_j},
\qquad
\phi_{L,j}^-(z)=\bigl(1-\frac{1}{z_j}\bigr)\prod_{l\not=j}\frac{z_j-tz_l}{z_j-z_l}
\end{split}
\end{equation*}
and with eigenvalue
\[E_\lambda^{qJ}=\sum_{j=1}^n\bigl(qabt^{2n-j-1}(q^{\lambda_j}-1)+t^{j-1}(q^{-\lambda_j}-1)\bigr).
\]
An algebraic formulation of the
orthogonality relations of the multivariable little $q$-Jacobi
polynomials is as follows. Define an $\Fz$-linear functional
\[h_L=h_L^{a,b;q,t}: \Fz[z]^S\rightarrow \Fz
\]
by requiring
\begin{equation*}
h_L\bigl(P_\lambda^L\bigr)=
\begin{cases}
1,\qquad &\hbox{ if } \lambda=0,\\
0,\qquad &\hbox{ if } \lambda\not=0.
\end{cases}
\end{equation*}
Define a form $\langle \cdot,\cdot\rangle_L= \langle
\cdot,\cdot\rangle_L^{a,b;q,t}$ on $\Fz[z]^S$ by
\begin{equation}\label{Lpairing}
\langle
p_1,p_2\rangle_L=h_L\bigl(p_1(z)\overline{p_2(z)}\bigr),\qquad
p_1,p_2\in \Fz[z]^S,
\end{equation}
where $p(z) \mapsto \overline{p(z)}$ stands for the anti-linear
algebra involution on $\Fz[z]$ which extends complex conjugation on
$\mathbb{C}$ by requiring that the parameters $a,b,c,d,q,t$ and the
variables $z_j$ to be formally real, e.g. $\overline{a}=a$ and
$\overline{z_j}=z_j$. The orthogonality relations now take the form
\begin{equation}
\langle P_\lambda^L,P_\mu^L\rangle_L=N_L(\lambda)\delta_{\lambda,\mu},
\qquad \lambda,\mu\in\Lambda_n
\end{equation}
for suitable explicit quadratic norms
$N_L(\lambda)=N_L(\lambda;a,b;q,t)\in \Fz$, which are recalled in Subsection \ref{little}.

The expressions
\begin{equation}\label{Dgeneralized}
D_q(\lambda;a,b;t)=
\frac{\bigl(P_\lambda^L(0;q^{-1}a,q^{-1}b;q,t)\bigr)^2}{N_L(\lambda;q^{-1}a,q^{-1}b;q,t)}
\in \Fz, \qquad \lambda\in\Lambda_n
\end{equation}
play a crucial role in the sequel. Note that $D_q(\lambda;a,b;t)$
is the generalized dimension $d(\lambda)$ from the introduction
(see \eqref{tildeD}), up to a slight reparametrization of the four
auxiliary parameters $a,b,q$ and $t$ (which we justify in a
moment).

We show that specializations of the parameters $(a,b,t)$ in $D_q(\lambda)$
give the quantum dimensions of spherical irreducible
representations associated to quantum complex Grassmannians.
In addition we show that suitable degenerations in $(a,b,q,t)$
lead to the dimensions of spherical irreducible representations
associated to complex/real Grassmannians (in which case the four
parameters degenerate to one, but with different exponents), and
to dimensions of spherical irreducible representations associated to
Grassmannians over non-Archimedean local fields (in which case $q$ degenerates to zero while the
remaining three parameters $(a,b,t)$ are specialized to $q$-independent values which are explicitly
given in terms of the cardinality of the residue field). The latter case thus
does not involve degenerating the parameters $(a,b,t)$, which is
the reason to use the multivariable little $q$-Jacobi polynomials with parameters $(q^{-1}a, q^{-1}b,t)$
in the definition \eqref{Dgeneralized} of the generalized dimensions.

We derive evaluation formulas for $P_\lambda^L$ in Section \ref{evaluation},
which thus lead to explicit evaluations of $D_q(\lambda)=D_q(\lambda;a,b;t)$
and to explicit evaluations of the (quantum) dimensions for the irreducible spherical
irreducible representations associated to different types of (quantum) Grassmannians.



\section{\textsf{CQG} algebras and quantum dimensions}\label{CQG}

The category of Compact Quantum Group algebras (\CQG) is a full
subcategory of the category of Hopf $*$-algebras that includes the
standard noncommutative deformations of algebras of representative
functions on compact Lie groups. In particular, \CQG algebras are
Hopf $*$-algebras with positive definite Haar functionals which have
natural analogs of the Peter-Weyl decomposition and the Schur
orthogonality relations. This leads to the notion of the quantum
dimension of an irreducible, finite dimensional comodule over a
given \CQG algebra. One of the main aims of the paper is to
explicitly compute quantum dimensions of irreducible comodules which
arise in the context of (quantum) Grassmannians. Following
\cite{DK}, we recall in this section the definition of a \CQG
algebra and its basic properties. We furthermore discuss the notion
of a quantum Gelfand pair in the context of \CQG algebras (cf. \cite{F}), and the associated
quantum zonal spherical functions.

Let $A$ be a Hopf $*$-algebra with comultiplication $\Delta$,
counit $\epsilon$ and antipode $S$. A linear functional
$h: A\rightarrow \mathbb{C}$ is called a normalized
Haar functional if $h(1)=1$ and
\[\bigl(h\otimes \hbox{Id}_A\bigr)\bigl(\Delta(a)\bigr)=
h(a)1_A=\bigl(\hbox{Id}_A\otimes h\bigr)\bigl(\Delta(a)\bigr),\qquad a\in A.
\]
A Haar functional $h$ on $A$ is called positive definite if
$h(a^*a)>0$ for all $0\not=a\in A$. There are different equivalent
definitions of Dijkhuizen's and Koornwinder's \cite{DK} notion of a
\CQG algebra, the following is a convenient starting point for our
purposes.
\begin{defi}
We call the Hopf $*$-algebra $A$ a Compact Quantum Group
\textup{(}\textsf{CQG}\textup{)} algebra if there exists a positive definite
Haar functional on $A$.
\end{defi}
The positive definite Haar functional on a \CQG algebra $A$ is
unique. Several concrete examples of \CQG algebras (such as the
algebra of representative functions on a compact group, the basic
commutative example) will be discussed in detail in the following
two sections.

Let $A$ be a \CQG algebra with normalized, positive
definite Haar functional $h$. The following results are all from
\cite{DK}. Any finite-dimensional right $A$-comodule is
unitarizable, hence completely reducible. We denote $\Sigma$ for the
set of equivalence classes of finite dimensional, (unitary)
irreducible right $A$-comodules. For $\alpha\in\Sigma$ let
$\pi_\alpha: L_\alpha\rightarrow L_\alpha\otimes A$ be a
representative of the corresponding irreducible, finite dimensional
$A$-comodule. Let $\pi_\alpha^*: L_\alpha^*\rightarrow
L_\alpha^*\otimes A$ be the corresponding dual right $A$-comodule,
defined by requiring
\[\bigl(\pi_\alpha^*(\phi)\bigr)\bigl(v\otimes 1_A)=
\bigl(\phi\otimes S\bigr)\bigl(\pi_\alpha(v)\bigr),\qquad v\in L_\alpha,\,\,\phi\in L_\alpha^*,
\]
where $1_A$ denotes the unit element of $A$.
Then $\{L_\alpha^*\otimes L_\alpha \, | \, \alpha\in\Sigma \}$
is a complete set of representatives of the irreducible, finite dimensional
$A\otimes A$-comodules.

Consider $A$ as a right $A\otimes A$-comodule by
\[R(a)=\sum_{(a)}a_{(2)}\otimes S(a_{(1)})\otimes a_{(3)},\qquad a\in A,
\]
where we use the standard Sweedler notation for (iterated)
applications of the comultiplication. The first key property of \CQG
algebras is the irreducible decomposition
\[A\simeq \bigoplus_{\alpha\in\Sigma}L_\alpha^*\otimes L_\alpha
\]
of $A$ as right $A\otimes A$-comodule (the Peter-Weyl Theorem for
\CQG algebras). Concretely, the isomorphism is realized by mapping
$\phi\otimes v\in L_\alpha^*\otimes L_\alpha$ to
$(\phi\otimes\hbox{Id}_A)(\pi_\alpha(v))\in A$.  We denote $A_h$ for
the Hilbert space completion of $A$ with respect to the pre-Hilbert
structure
\[\langle a,b\rangle_h:=h(b^*a).
\]

 Let $(\pi, V)$ be a finite dimensional, irreducible,
right $A$-comodule. Consider the canonical
linear isomorphism $\varphi: V\rightarrow V^{**}$ defined by
$\varphi(v)f:=f(v)$ for $f\in V^*$ and $v\in V$. We turn $V$ into a finite
dimensional Hilbert space such that $(\pi, V)$ is a unitary right
$A$-comodule (such scalar product is unique up to a strictly positive
scalar multiple). Since the double-dual right $A$-comodule
$(\pi^{**},V^{**})$ is equivalent to $(\pi, V)$,
there exists a unique positive definite, linear isomorphism
$F=F_V: V\rightarrow V$ such that
\[\hbox{Tr}_V(F)=\hbox{Tr}_V(F^{-1})>0
\]
and such that $\varphi\circ F: (\pi, V)\rightarrow (\pi^{**}, V^{**})$
is an intertwiner of right $A$-comodules (see \cite[Prop. 3.6]{DK}).

\begin{defi}\label{qDimdef}
Let $A$ be a \CQG algebra, $(\pi,V)$ a finite dimensional,
irreducible right $A$-comodule and $F: V\rightarrow V$ the linear
isomorphism as defined above. Then
\[\textup{Dim}_A(V):=\textup{Tr}_V\bigl(F\bigr)
\]
is called the quantum dimension of the right $A$-comodule $(\pi, V)$.
\end{defi}

A second key property of \CQG algebras is the generalized Schur
orthogonality relations \cite[Prop. 3.4 \& 3.5]{DK}, which we now
proceed to recall.

\begin{prop}\label{Schur}
Let $A$ be a \CQG algebra. Let $(\pi,V)$ and $(\pi',V')$ be finite
dimensional, irreducible, unitary right $A$-comodules with
associated scalar products $\langle \cdot,\cdot\rangle_{\pi}$ and
$\langle \cdot,\cdot\rangle_{\pi'}$.

\begin{enumerate}

\item [(a)] If $\pi \not\simeq\pi'$, then
\[\langle \langle \pi(u),v\rangle_{\pi},
\langle \pi'(u'), v'\rangle_{\pi'}\rangle_h=0
\]
for $u,v\in V$ and $u^\prime, v^\prime\in V'$.

\item [(b)] For $u_1,u_2,v_1,v_2\in V$ we have
\[\langle \langle \pi(u_1),v_1\rangle_\pi,
\langle \pi(u_2), v_2\rangle_\pi\rangle_h= \frac{\langle
u_1,u_2\rangle_\pi\langle F_V^{-1}v_2,v_1\rangle_\pi}
{\textup{Dim}_A(V)}.
\]
\end{enumerate}
\end{prop}

Let $A$ be a \CQG algebra. A quantum subgroup of $A$ is a pair
$(B,p)$ with $B$ a Hopf $*$-algebra and $p: A\rightarrow B$ a
surjective Hopf $*$-algebra homomorphism. In this situation, we
write
\[A_B=\{a\in A \, | \, \bigl(\hbox{Id}_A\otimes p\bigr)\bigl(\Delta(a)\bigr)=a\otimes 1_B \}
\]
for the $*$-subalgebra of right $B$-invariant elements in $A$. The
$*$-subalgebra $A_B\subset A$ is
a right $A$-comodule under the right regular co-action
\[\rho_B(a)=\sum_{(a)}a_{(2)}\otimes S(a_{(1)}),\qquad a\in A_B.
\]
By the Peter-Weyl Theorem for $A$, the right $A$-comodule $A_B$ decomposes into
irreducibles as
\begin{equation}\label{PWB}
A_B\simeq \bigoplus_{\alpha\in\Sigma} L_\alpha^*\,{}^{ \textup{Dim}_{\mathbb{C}}(L_\alpha^B)},
\end{equation}
where for a right $A$-comodule $(\pi, V)$,
\[V^B=\{v\in V \, | \, \bigl(1\otimes p)\bigl(\pi(v)\bigr)=v\otimes 1_B\}
\]
is the space of $B$-invariant elements in $V$.

\begin{defi}\label{qGelfandPair}
A \CQG algebra $A$ with quantum subgroup $(B,p)$ is called a quantum Gelfand
pair if
\[\textup{Dim}_{\mathbb{C}}(L_\alpha^B)\leq 1\qquad
\forall\,\alpha\in\Sigma,
\]
\textup{(}compare with \cite{F}\textup{)}.
\end{defi}
We fix a quantum Gelfand pair $(A,B)$. We call a finite
dimensional, irreducible $A$-comodule $B$-spherical if it has
nonzero $B$-invariant vectors. We denote $\Sigma_B\subset\Sigma$ for
the corresponding subset of equivalence classes of the $B$-spherical, irreducible,
finite dimensional right $A$-comodules.

In the various examples associated to (quantum) Grassmannians which
we discuss in the following sections, we will identify
$\Sigma_B$ with a certain fixed set of partitions. An important, but
case-by-case different, aspect is the concrete characterization of
the spherical irreducibles in terms of their parameterizing set
$\Sigma_B$. To maintain uniform notations as much as possible, it
will therefore be convenient to denote irreducible comodules labeled by $\Sigma_B$ by new symbols
$\{V_{\alpha}~|~ \alpha \in \Sigma_B \}$.

We write $\langle \cdot,\cdot\rangle_\alpha$ for a scalar product on
$V_\alpha$ that turns $V_\alpha$ into a unitary right $A$-comodule.
Note that $\Sigma_B$ is invariant under $\circ$, were $\circ$ is the
involution on $\Sigma$ such that $V_\alpha^*\simeq V_{\alpha^\circ}$
for all $\alpha\in \Sigma$. In fact, if we fix a nonzero vector
$v_\alpha\in V_\alpha^B$ ($\alpha\in\Sigma_B$), then $\langle
\cdot,v_\alpha\rangle_\alpha\in V_\alpha^*{}^B$. Combined with
\eqref{PWB}, we obtain the multiplicity free irreducible
decomposition
\begin{equation}\label{AB}
A_B=\bigoplus_{\alpha\in\Sigma_B}V_\alpha
\end{equation}
as right $A$-comodules.

{}From the present perspective, the natural definition of the
quantum zonal spherical functions would be
\[\widetilde{\varphi}_\alpha=\langle\pi_\alpha(v_\alpha), v_\alpha\rangle_\alpha, \qquad
\forall\,\alpha\in\Sigma_B,
\]
since they form a $\mathbb{C}$-basis of the $*$-algebra
\[{}_BA_B=\{a\in A \, | \, \bigl(\hbox{Id}_A\otimes p\otimes p\bigr)R(a)=a\otimes 1_B\otimes 1_B\}
\]
of $B$-biinvariant elements in $A$. It turns out though that the
following slightly modified definition of the quantum zonal spherical
functions is a better choice.

\begin{defi}\label{qzs}
The quantum zonal spherical functions
associated to the quantum Gelfand pair $(A,B)$ are defined by
\[\varphi_\alpha=\langle \pi_\alpha(v_\alpha), F_{V_{\alpha}}^{\frac{1}{2}}v_\alpha\rangle_\alpha\in A,\qquad
\forall\,\alpha\in\Sigma_B,
\]
 where $0\not=v_\alpha\in V_\alpha^B$ and $F_{V_\alpha}^{\frac{1}{2}}$ is the square
root of the positive definite linear operator $F_{V_\alpha}$ on $V_\alpha$.
\end{defi}
Note that the quantum zonal spherical function $\varphi_\alpha$
($\alpha\in\Sigma_B$) is uniquely determined up to nonzero scalar
multiples. The square root of the positive definite operator
$F_{V_\alpha}$ naturally appears in the definition of quantum zonal
spherical functions on quantum homogeneous spaces using Noumi's
\cite{N} two-sided coideal approach. It is this approach, and its
generalization in terms of coideal algebras by Letzter (see \cite{L}
and references therein), which has led to the identification of
Macdonald polynomials as radial parts of quantum zonal spherical
functions, see also Section \ref{quantum} for the quantum complex
Grassmannian.

It now follows from the preceding discussions that we have the
following expression for the quantum dimension of $V_\alpha$ in
terms of the associated quantum zonal spherical function
$\varphi_\alpha$ and the associated $B$-invariant vector $v_\alpha\in V_\alpha^B$.

\begin{prop}\label{propqdim}
Let $(A,B)$ be a quantum Gelfand pair. For
$\alpha\in\Sigma_B$ we have
\[\textup{Dim}_A\bigl(V_\alpha\bigr)=\frac{\|v_\alpha\|_\alpha^4}{\|\varphi_\alpha\|_h^2}.
\]
\end{prop}
\begin{proof}
This follows directly using the definition of $\varphi_\alpha$ and
the generalized Schur orthogonality relations (Proposition
\ref{Schur}{b}).
\end{proof}




\section{Grassmannians over local fields}\label{local}

In this section we specialize the general theory of the previous section to
the classical examples arising from Grassmannians over local fields.




\subsection{The Grassmann representation over local fields}\label{Gr}
Consider the commutative \CQG algebra $A=\RR(\Ko)$ of representative
functions on a compact topological group $\Ko$. The Hopf $*$-algebra
structures are given by
\[\Delta(f)(g,g^\prime)=f(gg^\prime), \quad \epsilon(f)=f(e),
\quad S(f)(g)=f(g^{-1}), \quad f^*(g)=\overline{f(g)},
\]
where $g,g^\prime\in \Ko$, $f\in \RR(\Ko)$, $e\in \Ko$ is the unit
element of $\Ko$, and where we have used the natural identification
$\RR(\Ko)\otimes \RR(\Ko)\simeq \RR(\Ko \times \Ko)$ by the
multiplication map. The positive definite normalized Haar functional
$h$ is
\[h(f)=\int_{\Ko} f(g)dg,\qquad f\in \RR(\Ko),
\]
where $dg$ is the normalized Haar measure on $\Ko$. Observe that
$\RR(\Ko)_h=L^2(\Ko,dg)$ and $\Sigma\simeq \widehat{\Ko}$, the
unitary dual of the compact group $\Ko$, since there is a natural
one-to-one correspondence between finite-dimensional right
$\RR(\Ko)$-comodules and finite dimensional $\Ko$-representations
preserving notions as irreducibility, equivalence and unitarity.
Concretely, for a given right $\RR(\Ko)$-comodule $\pi: V\rightarrow
V\otimes \RR(\Ko)$ we define a left $\Ko$-action on $V$ by
\[gv=\bigl(\hbox{Id}_V\otimes \hbox{ev}_g\bigr)\pi(v),\qquad g\in \Ko,\,\, v\in V,
\]
where $\hbox{ev}_g: \RR(\Ko)\rightarrow \mathbb{C}$ is the
evaluation map $\hbox{ev}_g(f)=f(g)$. Under this correspondence, the
dual representation of a finite dimensional
$\Ko$-representation $(\pi, V)$ becomes
\[
(\pi^*(g)f)(v)=f(\pi(g^{-1})v)\qquad  (f\in V^*, v\in V, g\in \Ko).
\]
Clearly the isomorphism $\pi^{**}\simeq \pi$ for a finite
dimensional $\Ko$-representation $(\pi, V)$ is realized by
$F=\hbox{Id}$, hence the quantum dimension
$\hbox{Dim}_{\RR(\Ko)}(V)$ is the complex dimension
$\hbox{Dim}_{\mathbb{C}}(V)$ of the representation space $V$. We now
define one of our main object of study.
\begin{defi}
Let $n,d\in\mathbb{N}$ with $n\leq d/2$. The Grassmannian
$\Gr(n,d;\Fb)$ over the local field $\Fb$ is the set of
$n$-dimensional $\Fb$-subspaces in $\Fb^d$.
\end{defi}
Consider the maximal compact subgroup $\KF$ of $\hbox{GL}_d(\Fb)$,
see \eqref{KF}. Note that $\KF$ acts transitively on $\Gr(n,d;\Fb)$,
e.g. by the Iwasawa decomposition in $\hbox{GL}_d(\Fb)$. The
subgroup of $\KF$ stabilizing $\{0\}^{d-n}\times \Fb^n \in
\Gr(n,d;\Fb)$ is denoted by $\LF$, so that
\[
\Gr(n,d;\Fb) \simeq \KF / \LF
\]
as transitive $\KF$-space. The pair $(\KF,\LF)$ is a Gelfand pair in the usual sense of
topological groups. This is well known for the Archimedean fields
$\mathbb{R}$ and $\mathbb{C}$, see \cite{JC} and \cite{DS} and the
references therein. For non-Archimedean fields we refer to
\cite{BO2}. Consequently, $(A,B)=\bigl(\RR(\KF), \RR(\LF)\bigr)$ is
a quantum Gelfand pair in the sense of \CQG algebras, where the surjective
Hopf $*$-homomorphism $p: \RR(\KF)\rightarrow \RR(\LF)$ is the
canonical restriction map $p(f)=f|_{\LF}$. The parametrizing set
$\Sigma_{\RR(\LF)}$ of the finite dimensional, irreducible
$\RR(\LF)$-spherical right $\RR(\KF)$-comodules is in bijection with
$(\KF/\LF)\hat{~}$, the unitary spherical dual of $\KF$ with respect
to the subgroup $\LF$. The corresponding right $\RR(\KF)$-comodule
$\RR(\KF)_{\RR(\LF)}$, viewed as $\KF$-representation, now
identifies with the {\it Grassmann representation}
\[
\RR\bigl(\Gr(n,d;\Fb)\bigr):=\{f\in \RR(\KF) \,\,\, | \,\,\, f(gh)=f(g),\,\,\,
\forall\,g\in \KF,\,\forall\, h\in \LF \},
\]
considered as left $\KF$-module by the regular action
\[\bigl(gf\bigr)(g^\prime)=f(g^{-1}g^\prime),\qquad f\in \RR\bigl(\Gr(n,d;\Fb)\bigr),\,\, g, g^\prime
\in \KF.
\]

In this setting we take $(\pi_\alpha^{\mathbb{F}}, V_\alpha^{\Fb})$
($\alpha\in(\KF/\LF)\hat{~}$) to be a complete set of
representatives of the irreducible
$\LF$-spherical $\KF$-representations. Then
\begin{equation}\label{abstractdecomposition}
\RR\bigl(\Gr(n,d;\Fb)\bigr)\simeq \bigoplus_{\alpha\in
(\KF/\LF)\hat{~}}V_\alpha^{\Fb}
\end{equation}
is the multiplicity-free decomposition of
$\RR\bigl(\Gr(n,d;\Fb)\bigr)$ in irreducible $\KF$-representations,
cf. \eqref{AB}. For $\alpha\in(\KF/\LF)\hat{~}$ we choose a scalar
product $\langle \cdot,\cdot\rangle_\alpha$ on
$V_\alpha^{\mathbb{F}}$ turning $\pi_\alpha^{\Fb}$ into a unitary
$\KF$-representation, and we choose a nonzero $\LF$-fixed vector
$v_\alpha^{\Fb}\in V_\alpha^{\Fb}$. Using again the identification
of finite dimensional right $\mathbb{C}(\KF)$-comodules with finite
dimensional $\KF$-representations, the quantum zonal spherical
functions associated to the quantum Gelfand pair $(\RR(\KF),
\RR(\LF))$ are the usual zonal spherical functions
\[\varphi_\alpha^{\Fb}(g)=\langle \pi_\alpha^{\Fb}(g)v_\alpha^{\Fb},
v_\alpha^{\Fb}\rangle_\alpha,\qquad g\in \KF,\,\,\alpha\in
(\KF/\LF)\hat{~}
\]
associated to the Gelfand pair $(\KF, \LF)$.
The dimension formula (see Proposition \ref{propqdim}) now
becomes
\begin{equation}\label{dimforclassical}
\hbox{Dim}_{\mathbb{C}}\bigl(V_\alpha^{\Fb}\bigr)=
\frac{\bigl(\varphi_\alpha^{\Fb}(e)\bigr)^2}{\int_{\KF}|\varphi_\alpha^{\mathbb{F}}(g)|^2dg},
\qquad \alpha\in (\KF/\LF)\hat{~},
\end{equation}
where $e$ is the unit element of $\KF$.



\subsection{Irreducible constituents and dimension formulas: the Archimedean
case}\label{cG}

The harmonic analysis corresponding to the complex and real
Grassmannian is classical by now. It relates to special cases of
harmonic analysis on compact symmetric spaces. To emphasize the
similarities with the non-Archimedean case and the quantum case, we
now shortly describe the relevant results for the complex
Grassmannian $\Gr(n,d;\mathbb{C})$, following closely the
presentation in \cite[Section 2]{DS}. For the real case, we refer
to \cite{JC} and especially to \cite{Su} (the presentation in \cite{Su} is very close
to the (quantum) complex case as discussed in this paper).

To avoid confusion with
standard notations for the complexification of a compact Lie group,
we write
\begin{equation}\label{convention}
(\Ko,\Le)=(\hbox{U}(d), \hbox{U}(d-n)\times\hbox{U}(n))
\end{equation}
for the Gelfand pair $(\KF, \LF)$ when $\Fb=\mathbb{C}$.
The irreducible $\Ko$-representations are parametrized by the cone
\begin{equation}\label{Pplus}
P_d^+=\{\mu\in\mathbb{Z}^{d} \, | \,
\mu_1\geq\mu_2\geq\cdots\geq\mu_{d}\}\subset \mathbb{Z}^d,
\end{equation}
by associating to $\mu\in P_d^+$ the $\Ko$-representation $L_\mu^\mathbb{C}$ of
highest weight $\mu\in P_d^+$. The Weyl character formula states
that the restriction of the character of $L_\mu^\mathbb{C}$ ($\mu\in
P_d^+$) to the compact $d$-torus of diagonal matrices
$\hbox{diag}(u_1,\ldots,u_d)$ in $\hbox{U}(d)$ is the Schur function
\[
s_\mu(u_1,\ldots,u_{d})=\Delta(u)^{-1}\sum_{w\in
S_{d}}(-1)^{l(w)}u^{w(\mu+\rho_d)}\quad (\mu\in P_d^+),
\]
where $\Delta(u)=\prod_{1\leq i<j\leq d}(u_i-u_j)$ is the Vandermonde determinant,
$l(w)$ is the length of $w\in S_n$, and $\rho_d=(d-1,d-2,\ldots,0)\in P_d^+$ is the staircase partition. The
dimension of $L_\mu^\mathbb{C}$ being the evaluation of its
character at the unit element $e\in \Ko$ thus yields
\begin{equation}\label{Schurdim}
\hbox{Dim}_{\mathbb{C}}(L_\mu^q)=s_\mu(1),\qquad \mu\in P_d^+.
\end{equation}
The resulting, famous, Weyl dimension formula is
\begin{equation}\label{Weyl}
\hbox{Dim}_{\mathbb{C}}(L_\mu^{\mathbb{C}})=\prod_{1\leq i<j\leq d}\left(
\frac{\mu_i-\mu_j+\delta_i-\delta_j}{\delta_i-\delta_j}\right),
\qquad \mu\in P_d^+,
\end{equation}
where
\begin{equation}\label{delta}
\delta=\frac{1}{2}\bigl(d-1,d-3,\ldots,3-d,1-d\bigr)
\end{equation}
is the half sum of positive roots in the standard realization
of the $A_{d-1}$-type root system in $\mathbb{R}^d$. We note that other
explicit dimension formulas are known, mostly with a combinatorial
flavour. For instance, for $\mu\in P_d^+$ with $\mu_d\geq 0$,
$\hbox{Dim}_{\mathbb{C}}(L_\mu^\mathbb{C})$ is the number of standard Young tableaux
of shape $\mu$. Furthermore, $\hbox{Dim}_{\mathbb{C}}(L_\mu^\mathbb{C})$
can be explicitly expressed in terms of products of hook-lengths of the Young diagram of shape $\mu$ (see \cite{M}).

We state now a different type of dimension formula for the subclass
of spherical irreducible representations associated to the complex
Grassmannian $\mathcal{G}(n,d;\mathbb{C})\simeq \Ko/\Le$. This
dimension formula is based on the expression \eqref{dimforclassical}
of the dimension in terms of the corresponding zonal spherical
functions. The spherical unitary dual $(\KF/\LF)\hat{~}$ naturally
identifies with the set $\Lambda_n$ of partitions of length at most
$n$ via the embedding
\begin{equation}\label{natural}
\lambda^\natural=
(\lambda_1,\lambda_2,\ldots,\lambda_n,\underset{d-2n}{\underbrace{0,\ldots,0}},-\lambda_n,\ldots,-\lambda_2,-\lambda_1)
\end{equation}
of $\Lambda_n$ in $P_d^+$. Consequently the irreducibles
\[V_\lambda^{\mathbb{C}}:=L_{\lambda^\natural}^{\mathbb{C}},\qquad
(\lambda\in\Lambda_n)
\]
form a complete set of representatives of the
irreducible $\LF$-spherical
$\KF$-re\-pre\-sen\-ta\-tions, and the Grassmann representation
$\RR\bigl(\Gr(n,d;\mathbb{C})\bigr)$ has a mul\-ti\-pli\-ci\-ty-free
decomposition
\[\RR\bigl(\Gr(n,d;\mathbb{C})\bigr)\simeq\bigoplus_{\lambda\in\Lambda_n}V_{\lambda}^{\mathbb{C}}
\]
as $\KF$-module. We write $\varphi_\lambda^{\mathbb{C}}$ for the
zonal spherical function corresponding to the irreducible
representation $V_{\lambda}^{\mathbb{C}}$ ($\lambda\in\Lambda_n$).
\begin{prop}\label{classicalcomplexformula}
For $\lambda\in\Lambda_n$ we have
\begin{equation*}
\begin{split}
\hbox{Dim}_{\mathbb{C}}\bigl(V_{\lambda}^{\mathbb{C}}\bigr)&=
\prod_{i=1}^n
\left(\frac{d-2n+1+2(\lambda_i+\rho_i)}{d-2n+1+2\rho_i}\right)\prod_{j=1}^{d-2n}\left(\frac{j+\lambda_i+\rho_i}{j+\rho_i}\right)^2\\
&\times\prod_{1\leq j<k\leq n}\left(\frac{d-2n+1+\lambda_j+\lambda_k+\rho_j+\rho_k}{d-2n+1+\rho_j+\rho_k}\right)^2
\left(\frac{\lambda_j-\lambda_k+\rho_j-\rho_k}{\rho_j-\rho_k}\right)^2,
\end{split}
\end{equation*}
where $\rho=(n-1,n-2,\ldots,1,0)\in\Lambda_n$ is the staircase partition.
\end{prop}
The natural proof of the proposition uses the identification of the radial parts of the zonal spherical
functions $\varphi_\lambda^{\mathbb{C}}$ ($\lambda\in\Lambda_n$) as $\hbox{BC}_n$-type Heckman-Opdam polynomials
with appropriate multiplicity parameters (see e.g. \cite{JC} and \cite[Section 2]{DS}).
The dimension formula \eqref{dimforclassical} then yields an expression in terms of evaluations and
quadratic norms of Heckman-Opdam polynomials, which both
have been explicitly evaluated.

Proposition \ref{classicalcomplexformula} follows also as a direct consequence of the formula
\begin{equation}\label{DclasQ}
\hbox{Dim}_{\mathbb{C}}\bigl(V_\lambda^{\mathbb{C}}\bigr)=
D_{q^2}\bigl(\lambda;q^{2(d-2n+1)},q^2;q^2\bigr)|_{q=1},\qquad
(\lambda\in\Lambda_n),
\end{equation}
expressing the dimensions in terms of the
classical ($q=1$) degeneration of the generalized dimension formula $D_q(\lambda)$ (see
\eqref{Dgeneralized}). Formula \eqref{DclasQ} is a direct consequence of the identification of multivariable little
$q$-Jacobi polynomials with quantum zonal spherical functions on the quantum complex Grassmannian,
see Subsection \ref{QDF} for details. We note also that the classical ($q=1$) degeneration of the
multivariable little $q$-Jacobi polynomials have been studied in
detail in \cite{SK}, also in the context of harmonic analysis on (quantum) complex Grassmannians (see \cite{DS}).

\begin{rema}
Proposition \ref{classicalcomplexformula}
can easily be reconfirmed by specializing the Weyl dimension formula \eqref{Weyl}
to $\mu=\lambda^\natural$ ($\lambda\in\Lambda_n$).
\end{rema}

For the real Grassmannian,
$(\Ko_{\mathbb{R}}/\Le_{\mathbb{R}})\hat{~}$ also naturally
identifies with $\Lambda_n$, see \cite{JC}, \cite{OU1} and \cite{Su}. The harmonic analytic results from
\cite{JC} lead to similar explicit dimension formulas for the associated
irreducible $\Le_{\mathbb{R}}$-spherical
$\Ko_{\mathbb{R}}$-representations $V_\lambda^{\mathbb{R}}$ ($\lambda\in\Lambda_n$). In particular
the dimensions can also be expressed as classical ($q=1$) degenerations
of generalized dimensions
\begin{equation}\label{realdim}
\textup{Dim}_{\mathbb{C}}\bigl(V_\lambda^{\mathbb{R}}\bigr)=
D_{q^2}\bigl(\lambda;q^{d-2n+1},q;q\bigr)|_{q=1},\qquad
(\lambda\in\Lambda_n),
\end{equation}
which follow as in the complex case from the results in \cite{JC} and \cite{SK} (see also \cite[\S 2.1.1]{OU1}).


\subsection{Irreducible constituents and dimension formulas: the non-Archimedean
case}\label{repp}

In this subsection we fix a non-Archimedean local field $\Fb$. Let
$\mathcal{O}$ denote the ring of integers and $\mathfrak{p}$ the
maximal ideal.
The Grassmann representation $\RR\bigl(\mathcal{G}(n,d;\Fb)\bigr)$ was studied
in \cite{Hill} and later in \cite{BO1,BO2}. Therein it is shown that
the irreducible $\LF$-spherical $\KF$-representations are also parameterized by
$\Lambda_n$. We now proceed to recall the explicit construction of
the irreducibles. The path to the construction of these
representations is completely different compared to the Archimedean
case. First of all, there is no (highest) weight theory, in
particular there is no concrete description of the full unitary dual
$\hat{\Ko}_{\Fb}$ as in the Archimedean case. We thus have to resort
to a direct construction of the spherical unitary dual
$(\KF/\LF)\hat{~}$ and the concrete construction of the
corresponding representations, which has a combinatorial flavor due
to the profinite nature of $\KF$,
\[\KF
\simeq\underset{\longleftarrow}{\textup{lim}}\,\textup{GL}_d(\mathcal{O}/\mathfrak{p}^k).
\]
The latter identification is realized by the canonical epimorphisms
$\pi_k: \KF \rightarrow \textup{GL}_d(\mathcal{O}/\mathfrak{p}^k)$
defined by reduction modulo $\mathfrak{p}^k$
($k\in\mathbb{Z}_{>0}$). Denoting $I_k$ for the kernel of $\pi_k$,
one has (cf. \cite{BO2})
\[
\RR\bigl(\Gr(n,d;\Fb)\bigr)\simeq
\underset{\longrightarrow}{\textup{lim}}\,\RR\bigl(\Gr(n,d;\Fb)\bigr)^{I_k}
\]
as $\KF$-representations, where $\RR\bigl(\Gr(n,d;\Fb)\bigr)^{I_k}$ stands for
the $\KF$-module consisting of the $I_k$-fixed vectors in the
Grassmann representation. The structure of the Grassmann
representation is completely determined by its so-called level $k-1$
components $\RR\bigl(\Gr(n,d;\Fb)\bigr)^{I_k}$ for $k\in\mathbb{Z}_{>0}$. We
proceed now to recall the precise structure of these level
components.

The set of partitions $\Lambda_n$ parametrizes the isomorphism
classes of the finite dimensional $\mathcal{O}$-modules of rank
$\leq n$. Concretely, the partition $\lambda\in\Lambda_n$
corresponds to the isomorphism class of the finite
$\mathcal{O}$-module $\oplus_{i=1}^n
\mathcal{O}/\mathfrak{p}^{\lambda_i}$. The finite
$\mathcal{O}$-modules in the isomorphism class labeled by
$\lambda\in\Lambda_n$ are called $\mathcal{O}$-modules of type
$\lambda$. Denote $\subseteq$ for the partial order on partitions
defined by inclusion of Young diagrams.  Fix a module of type
$k^d=(k,\ldots,k)$ ($d$ entries). For any type $\lambda \subseteq
k^d$ let $\Gr(\lambda,k^d;\mathcal{O})$ be the Grassmannian of its
submodules of type $\lambda$. Let
$\RR\bigl(\Gr(\lambda,k^d;\mathcal{O})\bigr)$ be the vector space of
complex valued functions on $\Gr(\lambda,k^d;\mathcal{O})$. The
canonical action of $\KF$ on $\Gr(\lambda,k^d;\mathcal{O})$ gives
rise to a representation of $\KF$ on
$\RR\bigl(\Gr(\lambda,k^d;\mathcal{O})\bigr)$. The finite Grassmann
representations relate to the Grassmann representation
$\RR\bigl(\Gr(n,d;\Fb)\bigr)$ by the identification
\begin{equation}\label{identification}
\RR\bigl(\Gr(n,d;\Fb)\bigr)^{I_k}\simeq\RR\bigl(\Gr(k^n,k^d;\mathcal{O})\bigr)
\end{equation}
as $\KF$-modules, see \cite[Lemma 2.1]{BO2}.

The irreducible constituents of the finite Grassmann representation
$\RR\bigl(\Gr(k^n,k^d;\mathcal{O})\bigr)$ have been explicitly described in
\cite{BO2} in terms of the so-called cellular basis of the
associated finite Hecke algebra
$\hbox{End}_{\KF}\bigl(\RR\bigl(\Gr(k^n,k^d;\mathcal{O})\bigr)\bigr)$. It
leads to a natural parametrization of the irreducible
$\LF$-spherical $\KF$-representations by partitions $\Lambda_n$, as
described in the following theorem, which is a direct consequence of
\cite[Theorem 1]{BO2}.

Denote $\Lambda_n^k$ for
partitions $\lambda\in\Lambda_n$ satisfying $\lambda\subseteq k^n$. Furthermore,
for an irreducible $\KF$-module $V$ we denote $\langle V:\RR\bigl(\Gr(\mu,k^d;\mathcal{O})\bigr)\rangle$
for the multiplicity of $V$ in the $\KF$-module $\RR\bigl(\Gr(\mu,k^d;\mathcal{O})\bigr)$.
\begin{thm}[\cite{BO2}]\label{pthm}
For $\lambda\in\Lambda_n$ there exists a unique irreducible
$\LF$-spherical $\KF$-re\-pre\-sen\-ta\-tion $V_\lambda^{\mathbb{F}}$
satisfying,
for all $k\in\mathbb{Z}_{>0}$ such that $\lambda\in\Lambda_n^k$,
\begin{enumerate}
\item[(1)] $\langle V_\lambda^{\mathbb{F}}:\RR\bigl(\Gr(\mu,k^d;\mathcal{O})\bigr)\rangle\geq 1$
for $\mu\in\Lambda_n^k$ satisfying $\lambda\subseteq\mu$,
\item[(2)] $\langle V_\lambda^{\mathbb{F}}:\RR\bigl(\Gr(\mu,k^d;\mathcal{O})\bigr)\rangle=0$
for $\mu\in\Lambda_n^k$ satisfying $\lambda\not\subseteq\mu$.
\end{enumerate}
The $\KF$-representations $V_\lambda^{\mathbb{F}}$ \textup{(}$\lambda\in\Lambda_n$\textup{)} form a complete
set of representatives of the irreducible $\LF$-spherical $\KF$-representations.
\end{thm}

Theorem \ref{pthm} gives
the multiplicity free irreducible decomposition
\[
\RR\bigl(\Gr(n,d;\Fb)\bigr)\simeq\bigoplus_{\lambda \in
\Lambda_n}V^{\Fb}_{\lambda}
\]
as $\KF$-modules, as well as the identifications
\[
(\KF/\LF)\hat{~}\simeq\Lambda_n\simeq\left\{ \text{isomorphism types of
$\mathcal{O}$-modules of rank $\le n$}  \right\}.
\]
Furthermore, \eqref{identification} and Theorem \ref{pthm} imply
\begin{equation}\label{irrk}
\RR\bigl(\Gr(k^n, k^d;\mathcal{O})\bigr)\simeq\bigoplus_{\lambda\in\Lambda_n^k}V^{\Fb}_{\lambda}
\end{equation}
for $k\in\mathbb{Z}_{>0}$ as $\KF$-modules, cf. \cite[Prop. 2.2]{BO2}.
The following is one of the main results in this paper.
\begin{thm}\label{thm-dim-padic}
 Let $t=|\mathcal{O}/\mathfrak{p}|^{-1}$. We have
\begin{equation}\label{padic-dimension}
\textup{Dim}_{\mathbb{C}}(V^{\Fb}_{\lambda})=t^{-(d-2n+1)|\lambda|-2(\rho,\lambda)}\left[
{n \atop \partial\lambda' }\right]_t
\frac{(t^{d-\lambda'_1-\lambda'_2+2};t)_{\lambda'_1+\lambda'_2}}{(t^{n-\lambda'_1+1};t)_{\lambda'_1}}
\frac{(1-t^{d-2\lambda'_1+1})}{(1-t^{d+1})}
\end{equation}
for $\lambda\in\Lambda_n$.
\end{thm}
Theorem \ref{thm-dim-padic} is a direct consequence of the formula
\begin{equation}\label{reldimD0}
\hbox{Dim}_{\mathbb{C}}\bigl(V_\lambda^{\mathbb{F}}\bigr)=
D_0\bigl(\lambda;t^{d-2n+1},t;t),\qquad
t=|\mathcal{O}/\mathfrak{p}|^{-1}
\end{equation}
from \cite{OU1}, and the explicit evaluation of $D_0(\lambda;a,b;t)$ derived in Theorem \ref{padicformula}.
Formula \eqref{reldimD0} is obtained in \cite{OU1} as
a consequence of \eqref{dimforclassical} and the identification of
the zonal spherical function $\varphi_\lambda^{\mathbb{F}}$ with a
suitable $p$-adic ($q=0$) degeneration of the multivariable little
$q$-Jacobi polynomial $P_\lambda^L$.

\begin{cor}\label{identity} For $k\in\mathbb{Z}_{>0}$ we have
\[
\sum_{\lambda \subseteq k^n}
t^{-(d-2n+1)|\lambda|-2(\rho,\lambda)}\left[ {n \atop
\partial\lambda' }\right]_t
\frac{(t^{d-\lambda'_1-\lambda'_2+2};t)_{\lambda'_1+\lambda'_2}}{(t^{n-\lambda'_1+1};t)_{\lambda'_1}}
\frac{(1-t^{d-2\lambda'_1+1})}{(1-t^{d+1})} =t^{-n(d-n)(k-1)}\left(
{d \atop n }\right)_{t^{-1}}.
\]
\end{cor}
\begin{proof}
Since the desired expression is rational in $t$, it suffices to
prove the identity with the formal parameter $t^{-1}$ specialized to
the cardinality of the residue field of a non-Archimedean local
field $\mathbb{F}$. In this case, the identity immediately follows
from \eqref{irrk}, \eqref{identification}, \eqref{padic-dimension}
and the easily verified identity (alternatively see e.g. \cite[Claim
4.1(1)]{BO2})
\begin{equation}\label{cardinality}
|\Gr(k^n,k^d;\mathcal{O})|=t^{-n(d-n)(k-1)}\left( {d \atop n
}\right)_{t^{-1}},\qquad t=|\mathcal{O}/\mathfrak{p}|^{-1}.
\end{equation}
\end{proof}
For $n=1$, Corollary \ref{identity} reduces to the trivial geometric
sum identity
\[
1+t^{1-d}\frac{(1-t^{d-1})}{(1-t)}+\frac{(1-t^{d-1})(1-t^d)}{(1-t)}
\sum_{m=2}^kt^{(1-d)m}=t^{(1-d)k}\frac{(1-t^d)}{(1-t)}.
\]

We end this subsection by discussing the dimension formula
\eqref{padic-dimension} in two cases that the dimensions have been
computed before by representation theoretic methods. The dimension
formula \eqref{padic-dimension} yields for $r=1,\ldots,n$,
\begin{equation}\label{padicfundamental}
\hbox{Dim}_{\mathbb{C}}\bigl(V_{\omega_r}^{\mathbb{F}}\bigr)
=t^{r(r-d)}\frac{\bigl(t^{d-r+2};t\bigr)_r(1-t^{d-2r+1})}
{\bigl(t;t\bigr)_r(1-t^{d+1})}=\left( {d \atop r}\right)_{t^{-1}}-\left({d \atop r-1}
\right)_{t^{-1}}
\end{equation}
where $t=|\mathcal{O}/\mathfrak{p}|^{-1}$.
Alternatively, by \eqref{irrk}
\[\RR\bigl(\mathcal{G}(\omega_r,1^d;\mathcal{O})\bigr)\simeq
\RR\bigl(\mathcal{G}(\omega_{r-1},1^d;\mathcal{O})\bigr)\oplus
V_{\omega_r}^{\Fb},
\]
which implies \eqref{padicfundamental} in view of \eqref{cardinality}. The latter
representation theoretic derivation essentially is the argument
from \cite{D}.

A similar analysis applies for the irreducible representations
$V_{k}^{\Fb}$ ($k\in \mathbb{Z}_{\geq 2}\subset \Lambda_1$) arising in the
representation $\RR\bigl(\Gr(1,d;\Fb)\bigr)$ associated to the projective
space $\mathbb{P}^d(\Fb)=\Gr(1,d;\Fb)$, see \cite{Haran} and
\cite[\S4]{OU1}. In this case, \eqref{padic-dimension} yields
\[
\hbox{Dim}_{\mathbb{C}}\bigl(V_{k}^{\Fb}\bigr)=t^{-(d-1)k}\frac{(1-t^{d-1})(1-t^d)}{(1-t)}
=t^{-(d-1)k}\frac{(1-t^d)}{(1-t)}-t^{-(d-1)(k-1)}\frac{(1-t^d)}{(1-t)}
\]
where $t=|\mathcal{O}/\mathfrak{p}|^{-1}$. This dimension formula can be reconfirmed using
\eqref{cardinality} and the decomposition
\[\RR\bigl(\mathcal{G}(k^1,k^d;\mathcal{O})\bigr)\simeq
\RR\bigl(\mathcal{G}((k-1)^1,(k-1)^d;\mathcal{O})\bigr)\oplus V_{k}^{\Fb}
\]
as $\KF$-modules, which is a direct consequence of \eqref{irrk}.


\section{The standard quantum complex Grassmannian}\label{qGrassmannian}


\subsection{The \CQG algebra of quantum representative functions on $\Ko$}

Let $\mathcal{U}_q$ be the quantized universal enveloping algebra of
$\mathfrak{gl}_n(\mathbb{C})$, where $0<q<1$. As a unital algebra it
is generated by the analogs $\xx_i,\yy_i$ ($i=1,\ldots,d-1$) of the
standard positive and negative simple root vectors of
$\mathfrak{gl}_n(\mathbb{C})$ and by the Cartan type elements
$\kk_j^{\pm 1}$ ($j=1,\ldots,d$). For the defining relations, as
well as for the explicit formulas turning $\mathcal{U}_q$ into a
Hopf algebra, we refer to \cite{OS}. We view $\mathcal{U}_q$ as a
Hopf $*$-algebra, with $*$-structure determined by
\[\xx_j^*=q^{-1}\yy_j\kk_j\kk_{j+1}^{-1},\qquad \yy_j^*=q\kk_j^{-1}\kk_{j+1}\xx_j,\qquad (\kk_i^{\pm 1})^*=\kk_i^{\pm 1}.
\]

A finite dimensional $\mathcal{U}_q$-representation decomposes in
common eigenspaces for the action of the Cartan type generators
$\kk_j$ ($j=1,\ldots,d$). We call a finite dimensional
$\mathcal{U}_q$-representation $V$ of type one if
\[V=\bigoplus_{\nu\in\mathbb{Z}^d}V[\nu],\qquad
V[\nu]=\{v\in V \,\, | \,\, \kk_jv=q^{\nu_j}v \quad (j=1,\ldots,d)
\}.
\]
We say that a vector $0\not=v\in V[\nu]$ has weight $\nu$. We define
$\RR_q(\Ko)$ to be the Hopf $*$-subalgebra of the Hopf dual of
$\mathcal{U}_q$ spanned by the matrix coefficients of the finite
dimensional $\mathcal{U}_q$-representations of type one. The Hopf
$*$-algebra structure on $\RR_q(\Ko)$ is obtained by dualizing the
Hopf $*$-algebra structure of $\mathcal{U}_q$, see \cite{OS}. In
particular, $\RR_q(\Ko)$ inherits a $*$-structure from $\mathcal{U}_q$
by
\[f^*(X)=\overline{f(S(X)^*)},\qquad f\in \RR_q(\Ko),\,\,\,
X\in\mathcal{U}_q,
\]
where $S$ is the antipode of $\mathcal{U}_q$. It is this particular
choice of $*$-structure on $\RR_q(\Ko)$ that reflects the fact that we
are considering quantized functions on the compact real form $\Ko$ of
$\hbox{GL}_d(\mathbb{C})$.

As in the classical case, which we discussed in Subsection \ref{cG}, there is
a one-to-one correspondence between finite dimensional (irreducible)
right $\RR_q(\Ko)$-comodules and finite dimensional (irreducible)
$\mathcal{U}_q$-representations of type one (see \cite{NYM}).
Concretely, for a finite dimensional right $\RR_q(\Ko)$-comodule
$\pi: V\rightarrow V\otimes \RR_q(\Ko)$, we define the associated
$\mathcal{U}_q$-action on $V$ by
\[\zz v=\bigl(\hbox{Id}\otimes \hbox{ev}_{\zz}\bigr)(\pi(v)),\qquad \zz \in\mathcal{U}_q,\,v\in V
\]
where $\hbox{ev}_{\zz}: \RR_q(\Ko)\rightarrow\mathbb{C}$ is the
evaluation map $\hbox{ev}_{\zz}(f)=f(\zz)$. Under this
identification, the $\mathcal{U}_q$-representation $V^*$ dual to $V$
is given by
\[
\bigl(\zz f\bigr)(v)=f\bigl(S(\zz)v),\qquad  \zz \in\mathcal{U}_q,\,
f\in V^*,\, v\in V.
\]
The equivalence classes of the finite dimensional, irreducible
$\mathcal{U}_q$-representations of type one are again naturally
parameterized by the cone $P_d^+$ \eqref{Pplus}. The corresponding
irreducible $\mathcal{U}_q$-representation $L^q_\mu$ ($\mu\in
P_d^+$) is characterized as the irreducible representation with
highest weight $\mu$. In particular, $L_\mu^q$ ($\mu\in P_d^+$) has
a (unique up to scalar multiples) weight vector $0\not=u_\mu\in
L_\mu^q[\mu]$ satisfying
\[
\xx_iu_\mu=0,\qquad i=1,\ldots,d-1,
\]
called a highest weight vector. Under the above mentioned
identification of right $\RR_q(\Ko)$-comodules and left
$\mathcal{U}_q$-modules, we have the Peter-Weyl decomposition
\begin{equation}\label{PWq}
\RR_q(\Ko)\simeq\bigoplus_{\mu\in P_d^+}L_\mu^q{}^*\otimes L_\mu^q,
\end{equation}
where $\phi\otimes v\in L_\mu^q{}^*\otimes L_\mu^q$ is identified
with the matrix coefficient $\phi(\cdot\,v)\in \RR_q(\Ko)$, which
gives the multiplicity free irreducible decomposition as right
$\RR_q(\Ko)\otimes \RR_q(\Ko)$-comodules (see \cite{NYM}). Consequently
$\RR_q(\Ko)$ is a \CQG algebra, with the normalized, positive definite
Haar functional $h$ on $\RR_q(\Ko)$ defined by requiring that it
vanishes on the matrix coefficients of the nontrivial irreducible
representations $L^q_\mu$ ($0\not=\mu\in P_d^+$). It follows from
the Peter-Weyl decomposition that $\Sigma\simeq P_d^+$, as for the
classical algebra of representative functions on $\Ko$. The quantum
Schur orthogonality relations (Proposition \ref{Schur}) for the \CQG
algebra $A=\RR_q(\Ko)$ are well known, see e.g. \cite[Section
3.2]{NYM}. To make the link precise, note that the quantum dimension
\[
\hbox{Dim}_q(V):=\hbox{Dim}_{\RR_q(\Ko)}(V)
\]
of a finite dimensional irreducible $\mathcal{U}_q$-representation
$V$ of type one has the following explicit form.

\begin{lem}\label{qdimalt}
Let $(\pi, V)$ be a finite dimensional irreducible
$\mathcal{U}_q$-representation of type one. Then
\[
\textup{Dim}_{q}(V)=\textup{Tr}_V(\kk^{2\delta}\cdot).
\]
\end{lem}
\begin{proof}
It is well known and easy to check that
\[S^2(\zz)=\kk^{2\delta}\zz \kk^{-2\delta},\qquad \forall\, \zz \in\mathcal{U}_q.
\]
Consequently, the isomorphism $\pi^{**}\simeq \pi$ as
$\mathcal{U}_q$-modules is realized by the linear operator
\[
F=\pi(\kk^{-2\delta})=\pi\bigl(\kk_1^{1-d}\kk_2^{3-d}\cdots
\kk_d^{d-1}\bigr).
\]
Turning $V$ into a $*$-unitary representation, the linear operator $F$ becomes
a positive definite linear operator on $V$ satisfying
$\hbox{Tr}_V(F)=\hbox{Tr}_V(F^{-1})$ (by an easy weight argument).
The result now follows from the definition of the quantum dimension
(Definition \ref{qDimdef}).
\end{proof}
\begin{rema}
The quantized universal enveloping algebra $\mathcal{U}_q$ is a
ribbon algebra. Lemma \ref{qdimalt} implies that the quantum
dimension $\hbox{Dim}_q(V)$ of an irreducible, finite dimensional
$\mathcal{U}_q$-module coincides with the natural (topologically
motivated) notion of quantum dimension for finite dimensional
modules over a ribbon algebra, see \cite[Chapter XIV]{Ka}.
\end{rema}
Since
$\hbox{Dim}_{\mathbb{C}}(L_\mu^q[\nu])=\hbox{Dim}_{\mathbb{C}}(L_\mu^{\mathbb{C}}[\nu])$
for all $\nu\in\mathbb{Z}^d$, one has
\begin{equation}\label{qdimSchur}
\hbox{Dim}_q(L_\mu^q)=\sum_{\nu\in\mathbb{Z}^d}\hbox{Dim}_{\mathbb{C}}(L_\mu^{\mathbb{C}}[\nu])q^{2\langle\delta,\nu\rangle}=
s_\mu(q^{2\delta}),\qquad \mu\in P_d^+.
\end{equation}
This implies the following $q$-analogue of the Weyl dimension
formula,
\begin{equation}\label{qWeyl}
\hbox{Dim}_q(L_\mu^q)=q^{-2\langle\delta,\mu\rangle}\prod_{1\leq
i<j\leq d}\left(\frac{1-q^{2(\mu_i-\mu_j+\delta_i-\delta_j)}}
{1-q^{2(\delta_i-\delta_j)}}\right),
\end{equation}
see e.g. \cite[Lemma 2.5]{St} for a simple proof based on the
classical Weyl denominator formula. A combinatorial formula for the
quantum dimension $\hbox{Dim}_q(L_\mu^q)$ ($\mu\in P_d^+$) is
provided by the $q$-analogue of the hook-length formula, see
\cite{M} and \cite[Section 3.2]{NYM} for more details.


\subsection{The quantum Grassmann representation}\label{qGr}

In this subsection we discuss the standard quantum analog of the
complex Grassmannian $\mathcal{G}(n,d;\mathbb{C})\simeq \Ko/\Le$ and
its associated quantum Grassmann representation $\RR_q(\Ko/\Le)$.
The one-parameter family of quantum complex Grassmannians from
\cite{NDS} is discussed separately in Section \ref{quantum}.

We define a Gelfand pair $(\RR_q(\Ko), \RR_q(\Le))$ of \CQG algebras
as follows. Write
$\mathfrak{m}^{\mathbb{C}}\simeq\mathfrak{gl}_{d-n}(\mathbb{C})\times\mathfrak{gl}_n(\mathbb{C})$
for the complexified Lie algebra of $\Le$. We view
$\mathfrak{m}^{\mathbb{C}}$ as Lie subalgebra of
$\mathfrak{gl}_d(\mathbb{C})$ in the standard way and write
\[\mathcal{U}_q(\mathfrak{m}^{\mathbb{C}})\simeq\mathcal{U}_q(\mathfrak{gl}_{d-n}(\mathbb{C}))\otimes\mathcal{U}_q(\mathfrak{gl}_n(\mathbb{C}))
\]
for the Hopf $*$-subalgebra of
$\mathcal{U}_q=\mathcal{U}_q(\mathfrak{gl}_d(\mathbb{C}))$ generated
by $\kk_i^{\pm 1}$ ($i=1,\ldots,d$) and $\xx_j,\yy_j$
($j\in\{1,\ldots,d-1\}\setminus \{d-n\}$). Let $\RR_q(\Le)$ be the
span of the matrix coefficients of the finite dimensional
irreducible
$\mathcal{U}_q(\mathfrak{m}^{\mathbb{C}})$-representations of type
one. Then $(\RR_q(\Le),p)$ is a quantum subgroup of $\RR_q(\Ko)$,
with $p$ the surjective Hopf $*$-algebra morphism
$p:\RR_q(\Ko)\rightarrow \RR_q(\Le)$ defined by
$p(f)=f|_{\mathcal{U}_q(\mathfrak{m}^{\mathbb{C}})}$.

Suppose $V$ is a finite dimensional right $\RR_q(\Ko)$-comodule.
Using the identification of right $\RR_q(\Ko)$-comodules and left
$\mathcal{U}_q$-modules, the subspace of $\RR_q(\Le)$-invariant
elements in $V$ coincides with the subspace of
$\mathcal{U}_q(\mathfrak{m}^{\mathbb{C}})$-invariant elements in
$V$,
\[V^{\RR_q(\Le)}=\{v\in V \,\, | \,\, \zz v=\epsilon(\zz)v\quad \forall \zz \in \mathcal{U}_q(\mathfrak{m}^{\mathbb{C}}) \}.
\]
In translating the general \CQG algebra theory from right
$\RR_q(\Ko)$-modules to left $\mathcal{U}_q$-re\-pre\-sen\-ta\-tions, we
shall refer to $\RR_q(\Ko)$-spherical representations as
$\mathcal{U}_q(\mathfrak{m}^{\mathbb{C}})$-spherical
representations. The branching rules for finite dimensional
$\mathcal{U}_q$-representations of type one, when viewed as
representations of $\mathcal{U}_q(\mathfrak{m}^{\mathbb{C}})\subset
\mathcal{U}_q$, are the same as for the classical Gelfand pair
$(\Ko,\Le)$. Consequently $(\RR_q(\Ko), \RR_q(\Le))$ is a quantum Gelfand pair,
with associated unitary spherical dual $\Sigma_{\RR_q(\Le)}\simeq (\Ko/\Le)\hat{~}
\simeq \Lambda_n$ (compare with the classical complex Grassmannian
of Subsection ref{cG}, e.g. \eqref{natural}). In particular, the
$\mathcal{U}_q$-representations
\[
V_\lambda^q:=L_{\lambda^\natural}^q,\quad (\lambda\in\Lambda_n)
\]
form a complete set of representatives of the finite dimensional,
irreducible $\mathcal{U}_q(\mathfrak{m}^{\mathbb{C}})$-spherical
$\mathcal{U}_q$-representations. The right $\RR_q(\Ko)$-comodule
$\RR_q(\Ko)_{\RR_q(\Le)}$ now identifies with the {\it quantum
Grassmann representation}
\[
\RR_q(\Ko/\Le):=\{f\in \RR_q(\Ko) \,\,\, | \,\,\, f(\zz
 Y)=f(\zz)\qquad \forall\, \zz \in\mathcal{U}_q, \forall\,
Y\in\mathcal{U}_q(\mathfrak{m}^{\mathbb{C}}) \},
\]
viewed as left $\mathcal{U}_q$-module by the regular
$\mathcal{U}_q$-action
\[\bigl( \zz f\bigr)(\zz ^\prime)=f(S(\zz) \zz^\prime),\qquad
f\in \RR_q(\Ko/\Le),\,\,\zz,\zz^\prime\in\mathcal{U}_q.
\]
It has the multiplicity free decomposition
\[
\RR_q(\Ko/\Le)\simeq\bigoplus_{\lambda\in\Lambda_n}V_\lambda^q
\]
in irreducible $\mathcal{U}_q$-representations.

For $\lambda\in\Lambda_n$ we fix a nonzero
$\mathcal{U}_q(\mathfrak{m}^{\mathbb{C}})$-invariant vector
$v_\lambda^q\in V_\lambda^q$, and we choose a scalar product
$\langle\cdot,\cdot\rangle_\lambda$ on $V_\lambda^q$ which turns
$V_\lambda^q$ into a $*$-unitary $\mathcal{U}_q$-representation. The
quantum zonal spherical functions $\varphi_\lambda^q\in
\RR_q(\Ko/\Le)$ ($\lambda\in\Lambda_n$) are given by
\[
\varphi_\lambda^q(\zz)=\langle \zz v_\lambda^q,
\kk^{-\delta}v_\lambda^q\rangle_\lambda,\qquad \zz\in\mathcal{U}_q,
\]
where $\kk^{\pm\delta}$ is the linear endomorphism of $V_\lambda^q$
satisfying
\[
\kk^{\pm\delta}|_{V_\lambda^q[\nu]}=q^{\pm\langle
\delta,\nu\rangle}\hbox{Id}_{V_\lambda^q[\nu]},\qquad \forall
\nu\in\mathbb{Z}^d
\]
(which is the square root of the positive definite linear operator
on $V_\lambda^q$ defined by the action of $\kk^{\pm 2\delta}$ on
$V_\lambda^q$). We define for
$\psi\in\hbox{End}_{\mathbb{C}}(V_\lambda^q)$,
\[
\varphi_\lambda^q(\psi):=\langle \psi(v_\lambda^q),
K^{-\delta}v_\lambda^q\rangle_\lambda,
\]
in particular
$\varphi_\lambda^q(\kk^{\delta})=\|v_\lambda^q\|_\lambda^2$. The
most convenient form of the quantum dimension formula (Proposition
\ref{propqdim}) then reads
\begin{equation}\label{dimforquantum}
\hbox{Dim}_q\bigl(V_\lambda^q\bigr)=\frac{\bigl(\varphi_\lambda^q(\kk^{\delta})\bigr)^2}{
\| \varphi_\lambda^q\|_h^2}, \qquad \lambda\in\Lambda_n.
\end{equation}


\subsection{Quantum dimension formulas}\label{QDF}

The quantum dimensions of the irreducible spherical $\mathcal{U}_q$-representations
$V_\lambda^q$ ($\lambda\in\Lambda_n$) have the following explicit closed expressions.
\begin{prop}\label{quantumcomplexformula}
For $\lambda\in\Lambda_n$ we have
\begin{equation}\label{qcf}
\begin{split}
\hbox{Dim}_{q}\bigl(V_{\lambda}^{q}\bigr)&=q^{2(2n-d-1)|\lambda|-4(\rho,\lambda)}
\prod_{i=1}^n\frac{\bigl(q^{2(1+\lambda_i+\rho_i)};q^2\bigr)_{d-2n}^2}
{\bigl(q^{2(1+\rho_i)};q^2\bigr)_{d-2n}^2}
\frac{(1-q^{2(d-2n+1+2(\lambda_i+\rho_i))})}
{(1-q^{2(d-2n+1+2\rho_i)})}\\
&\qquad\times\prod_{1\leq j<k\leq n}
\frac{(1-q^{2(d-2n+1+\lambda_j+\lambda_k+\rho_j+\rho_k)})^2
(1-q^{2(\lambda_j-\lambda_k+\rho_j-\rho_k)})^2}
{(1-q^{2(d-2n+1+\rho_j+\rho_k)})^2(1-q^{2(\rho_j-\rho_k)})^2}.
\end{split}
\end{equation}
\end{prop}
In Section \ref{quantum} we derive Proposition \ref{quantumcomplexformula} from a generalization of the quantum dimension formula
\eqref{dimforquantum} to the setup of the one-parameter family of quantum complex Grassmannians from \cite{NDS} and \cite{DS}.
By a suitable degeneration procedure we subsequently deduce the quantum dimension formulas
\begin{equation}\label{qcfcons}
\hbox{Dim}_q\bigl(V_\lambda^q\bigr)=
D_{q^2}\bigl(\lambda;q^{2(d-2n+1)},q^2;q^2),\qquad
\lambda\in\Lambda_n
\end{equation}
in terms of evaluations and quadratic
norms of multivariable little $q$-Jacobi polynomials. More directly, \eqref{qcfcons}
is a direct consequence of the quantum dimension formula \eqref{dimforquantum} if one uses the fact
that the quantum zonal spherical functions $\varphi_\lambda^q$ ($\lambda\in\Lambda_n$) identify with
multivariable little $q$-Jacobi polynomials (see \cite{DS} and Remark \ref{standardrema}).
\begin{rema}\label{classicalrema}
Taking the limit $q\rightarrow 1$ in \eqref{qcf}, we obtain the explicit dimension formulas
for the irreducible $\Le$-spherical $\Ko$-representations $V_\lambda^{\mathbb{C}}$
($\lambda\in\Lambda_n$) as stated in Proposition \ref{classicalcomplexformula}.
\end{rema}

We end this subsection by discussing the representation theoretic
viewpoint on the quantum dimension formula \eqref{qcf} for
$\lambda=\omega_r$, in which case it yields the explicit expression
\begin{equation}\label{dimformulanewrepr}
\hbox{Dim}_q\bigl(V_{\omega_r}^q\bigr)= q^{2r(d-r)}\left({d\atop
r}\right)_{q^{-2}}^2- q^{2(r-1)(d-r+1)}\left({d\atop
r-1}\right)_{q^{-2}}^2
\end{equation}
for $r=1,\ldots,n$ by a straightforward computation. The following
representation theoretic derivation of \eqref{dimformulanewrepr}
bears strong resemblance to the representation theoretic derivation
of the dimension formula \eqref{padicfundamental} in the classical
non-Archimedean setup.

Let $\{\epsilon_j\}_{j=1}^d$ be the standard orthonormal basis of
$\mathbb{R}^d$. Let $V=L_{\epsilon_1}^q$ be the $d$-dimensional
vector-representation of $\mathcal{U}_q$. Its weight spaces are
one-dimensional, with weights $\epsilon_j$ ($j=1,\ldots,d$). The
irreducible representation $L_{\epsilon_1+\cdots +\epsilon_r}^q$
($1\leq r\leq n$) can be realized as the $r$th graded part
$\Lambda_q^r(V)$ of the $q$-exterior algebra $\Lambda_q(V)$ of $V$
(see \cite{NYM}). The weight spaces of $\Lambda_q^r(V)$ are again
one-dimensional, with weights given by $\sum_{j\in J}\epsilon_j$ for
subsets $J\subseteq \{1,\ldots,d\}$ of cardinality $r$. Hence
\[\hbox{Dim}_q\bigl(\Lambda_q^r(V)\bigr)=
\hbox{Tr}_{\Lambda_q^r(V)}(\kk^{2\delta}\cdot)=
\sum_{\stackrel{J\subseteq\{1,\ldots,d\}}{\#J=r}}\prod_{j\in
J}q^{2(d-j)}=
q^{(r-1)r}\frac{\bigl(q^{2(d-r+1)};q^2\bigr)_r}{\bigl(q^2;q^2\bigr)_r},
\]
where the third equality follows by an easy induction argument. The
dual of $\Lambda_q^r(V)$ is irreducible of highest weight
$-\epsilon_{d-r+1}\cdots -\epsilon_d$. In fact, by \cite{NYM} we
have
\[
\Lambda_q^r(V)^*\simeq \mathbb{C}\hbox{det}_q^{-1}\otimes
\Lambda_{q}^{d-r}(V),
\]
with $\mathbb{C}\hbox{det}_q^{-1}\simeq L_{-\epsilon_1\cdots
-\epsilon_d}^q$ the one-dimensional representation realized by the
inverse of the quantum determinant $\det_q\in \RR_q(\Ko)$, hence
\[\hbox{Dim}_q\bigl(\Lambda_q^r(V)^*\bigr)=
q^{(1+r-2d)r}\frac{\bigl(q^{2(r+1)};q^2\bigr)_{d-r}}{\bigl(q^2;q^2\bigr)_{d-r}}.
\]
We conclude that
\[\hbox{Dim}_q\bigl(U_r\bigr)=
q^{2r(d-r)}\left({d\atop r}\right)_{q^{-2}}^2,\] where $U_r$ is the
$\mathcal{U}_q$-module
\[U_r=\Lambda_q^r(V)\otimes \Lambda_q^r(V)^*.
\]
By \cite[(6.14)]{DS} we have
\[U_r\simeq U_{r-1}\oplus V_{\omega_r}^q, \qquad r=1,\ldots,n
\]
as $\mathcal{U}_q$-modules (with $U_0$ the trivial representation),
which now immediately implies \eqref{dimformulanewrepr}.


\section{Evaluation formulas}\label{evaluation}

In this section we derive evaluation formulas for multivariable little and
big $q$-Jacobi polynomials by degenerating evaluation formulas for
the Macdonald-Koornwinder polynomials \cite{K}. Evaluation formulas for Macdonald-Koornwinder
polynomials have been obtained in \cite{vD}, \cite{Sa} and \cite{S2}.
The evaluation formulas for the multivariable little $q$-Jacobi
polynomials are used in the next section to obtain explicit expressions for the
generalized dimensions $D_q$ (see \eqref{Dgeneralized}) and their degenerations.


\subsection{Macdonald-Koornwinder polynomials}


Recall our notation $\Fz=\mathbb{C}(a,b,c,d,q,t)$ for the field of
rational functions in six indeterminates $a,b,c,d,q$ and $t$. Let
$W=S_n\ltimes \{\pm 1\}^n$, acting on $\mathbb{Z}^n$ by permutations
and sign-changes of the coordinates, and acting on the algebra
$\Fz[z^{\pm 1}]=\Fz[z_1^{\pm 1},\ldots,z_n^{\pm 1}]$ by permutations
and inversions of the independent variables $z_i$. The cone
$\Lambda_n\subset \mathbb{Z}^n$ is a fundamental domain for the
$W$-action on $\mathbb{Z}^n$. The dominance order on $\Lambda_n$ is
the partial order $\leq$ defined by $\lambda\leq\mu$ if
$\sum_{j=1}^{i}\lambda_j\leq\sum_{j=1}^i\mu_j$ for $1 \le i \le n$.
An $\Fz$-basis of the algebra $\Fz[z^{\pm 1}]^W$ of $W$-invariant
Laurent polynomials is given by the symmetric monomials
$m_\lambda(z):=\sum_{\mu\in W\lambda}z^\mu$ ($\lambda\in
\Lambda_n$).

The monic Macdonald-Koornwinder polynomial $P_\lambda(z)=P_\lambda(z;a,b,c,d;q,t)\in \Fz[z^{\pm 1}]^W$ of degree $\lambda\in\Lambda_n$
is of the form
\begin{equation}\label{expansionK}
P_\lambda(z)=m_\lambda(z)+\sum_{\mu\in\Lambda_n: \mu<\lambda}c_{\lambda,\mu}m_{\mu}(z)
\end{equation}
for certain coefficients $c_{\lambda,\mu}=c_{\lambda,\mu}(a,b,c,d;q,t)\in \Fz$. It can be characterized
as solution of a particular second-order difference equation, or, for suitably specialized generic parameters $a,b,c,d,q,t$,
in terms of suitable orthogonality properties, see Koornwinder \cite{K}.
The characterizing difference equation for the Macdonald-Koornwinder polynomial $P_\lambda(z)$ of degree $\lambda$ is
$DP_\lambda=E_\lambda P_\lambda$ where
\begin{equation*}
\begin{split}
D&=\sum_{j=1}^n\bigl(\phi_j(z)(T_j-\hbox{Id})+\phi_j(z^{-1})(T_j^{-1}-\hbox{Id})\bigr),\\
\phi_j(z)&=\frac{(1-az_j)(1-bz_j)(1-cz_j)(1-dz_j)}{(1-z_j^2)(1-qz_j^2)}
\prod_{l\not=j}\frac{(1-tz_lz_j)(1-tz_l^{-1}z_j)}{(1-z_lz_j)(1-z_l^{-1}z_j)}
\end{split}
\end{equation*}
and where the eigenvalue $E_\lambda\in \Fz$ is given by
\[E_\lambda=\sum_{j=1}^n\bigl(q^{-1}abcdt^{2n-j-1}(q^{\lambda_j}-1)+t^{j-1}(q^{-\lambda_j}-1)\bigr).
\]
Note that $P_\lambda(z;a,b,c,d;q,t)$ is symmetric in the four
parameters $a,b,c,d$. The evaluation formulas for
Macdonald-Koornwinder polynomials are given by
\begin{equation}\label{evalK}
\begin{split}
P_\lambda(at^\rho)&=\prod_{i=1}^n\frac{\bigl(abt^{n-i},act^{n-i},adt^{n-i},q^{-1}abcdt^{n-i};q\bigr)_{\lambda_i}}
{\bigl(q^{-1}abcdt^{2(n-i)};q\bigr)_{2\lambda_i}\bigl(at^{n-i}\bigr)^{\lambda_i}}\\
&\times\prod_{1\leq j<k\leq n}
\frac{\bigl(q^{-1}abcdt^{2n-j-k+1};q\bigr)_{\lambda_j+\lambda_k}
\bigl(t^{k-j+1};q\bigr)_{\lambda_j-\lambda_k}}
{\bigl(q^{-1}abcdt^{2n-j-k};q\bigr)_{\lambda_j+\lambda_k}
\bigl(t^{k-j};q\bigr)_{\lambda_j-\lambda_k}}
\end{split}
\end{equation}
for $\lambda\in\Lambda_n$. The evaluation formulas \eqref{evalK} have been established by van Diejen \cite{vD} for a
sub-family of Macdonald-Koornwinder polynomials. The general case follow from Sahi's \cite{Sa} results using double affine
Hecke algebras, see also \cite[Remark 9.5]{S2}.

In the rank one case ($n=1$), the Macdonald-Koornwinder polynomials do not dependent on $t$ and
reduce to the monic Askey-Wilson polynomials
\[P_m(z) =\frac{\bigl(ab,ac,ad;q\bigr)_m}{a^m\bigl(q^{m-1}abcd;q\bigr)_m}
{}_4\phi_3\left(\begin{matrix} q^{-m}, q^{m-1}abcd, az, az^{-1}\\ ab, ac, ad \end{matrix};q,q\right)\in\mathbb{C}[z+z^{-1}],
\qquad m\in\mathbb{Z}_+.
\]
In particular
\[P_m(a)=\frac{\bigl(ab,ac,ad;q\bigr)_m}{a^m\bigl(q^{m-1}abcd;q\bigr)_m},\qquad m\in\mathbb{Z}_+,
\]
which is in accordance to \eqref{evalK}.

We give here a purely algebraic formulation of the quadratic norm evaluations for the Macdonald-Koornwinder
polynomials, which has the advantage that we do not need to specialize the parameters $a,b,c,d,q$
and $t$. We define an $\Fz$-linear functional
\[h_K=h_K^{a,b,c,d;q,t}: \Fz[z^{\pm 1}]^W\rightarrow \Fz
\]
by requiring
\begin{equation*}
h_K\bigl(P_\lambda\bigr)=
\begin{cases}
1,\qquad &\hbox{ if } \lambda=0,\\
0,\qquad &\hbox{ if } \lambda\not=0.
\end{cases}
\end{equation*}
An analytic definition of $h_K$ can be given (for specific Zariski
dense choice of parameters) as an integral over a deformed compact
$n$-torus with explicit weight function, see e.g. \cite{K}.
We define now a sesqui-linear form $\langle \cdot,\cdot\rangle_K=
\langle \cdot,\cdot\rangle_K^{a,b,c,d;q,t}$
on $\Fz[z^{\pm 1}]^W$ by
\begin{equation}\label{Kpairing}
\langle
p_1,p_2\rangle_K=h_K\bigl(p_1(z)\overline{p_2(z)}\bigr),\qquad
p_1,p_2\in \Fz[z^{\pm 1}]^W,
\end{equation}
where complex conjugation on $\mathbb{C}$ is extended to an
anti-linear algebra involution on $\Fz[z^{\pm 1}]$ by requiring the
six parameters $a,b,c,d,q,t$ to be formally real, e.g.
$\overline{a}=a$, and by requiring the variables $z_j$ to be
formally purely imaginary, $\overline{z_j}=z_j^{-1}$. The
orthogonality relations for the Macdonald-Koornwinder polynomials
are now given by
\begin{equation}\label{normKoorn}
\langle P_\lambda,P_\mu\rangle_K=N_K(\lambda)\delta_{\lambda,\mu},\qquad
\forall\,\lambda,\mu\in\Lambda_n,
\end{equation}
with the quadratic norms $N_K(\lambda)=N_K(\lambda;a,b,c,d;q,t)$
given by $N_K(\lambda)=N_K^+(\lambda)N_K^-(\lambda)$, where
$N_K^{\pm}(\lambda)=N_K^{\pm}(\lambda;a,b,c,d;q,t)\in \Fz$ are
\begin{equation*}
\begin{split}
N_K^+(\lambda)&=\prod_{i=1}^n\frac{\bigl(abt^{n-i},act^{n-i},adt^{n-i},q^{-1}abcdt^{n-i};q\bigr)_{\lambda_i}}
{\bigl(q^{-1}abcdt^{2(n-i)};q\bigr)_{2\lambda_i}}\\
&\times\prod_{1\leq j<k\leq
n}\frac{\bigl(q^{-1}abcdt^{2n-j-k+1};q\bigr)_{\lambda_j+\lambda_k}
\bigl(t^{k-j+1};q\bigr)_{\lambda_j-\lambda_k}}
{\bigl(q^{-1}abcdt^{2n-j-k};q\bigr)_{\lambda_j+\lambda_k}\bigl(
t^{k-j};q\bigr)_{\lambda_j-\lambda_k}},\\
N_K^-(\lambda)&=\prod_{i=1}^n\frac{\bigl(qt^{n-i},bct^{n-i},bdt^{n-i},cdt^{n-i};q\bigr)_{\lambda_i}}
{\bigl(abcdt^{2(n-i)};q\bigr)_{2\lambda_i}}\\
&\times\prod_{1\leq j<k\leq
n}\frac{\bigl(abcdt^{2n-j-k-1};q\bigr)_{\lambda_j+\lambda_k}
\bigl(qt^{k-j-1};q\bigr)_{\lambda_j-\lambda_k}}
{\bigl(abcdt^{2n-j-k};q\bigr)_{\lambda_j+\lambda_k}\bigl(
qt^{k-j};q\bigr)_{\lambda_j-\lambda_k}}.
\end{split}
\end{equation*}
The quadratic norm formulas have been established by van Diejen \cite{vD} for a
sub-family of Macdonald-Koornwinder polynomials. The general case has been derived in \cite{Sa},
see also \cite{S2}.


\subsection{Multivariable little $q$-Jacobi polynomials}\label{little}

We now return to the multivariable little $q$-Jacobi
polynomials $P_\lambda^L(z)=P_\lambda^L(z;a,b;q,t)$ from Section \ref{qJsection}.
We first shortly discuss evaluation formulas for one-variable little $q$-Jacobi polynomials ($n=1$),
in which case the evaluation formulas can be derived directly from the explicit expressions of the little $q$-Jacobi polynomials
as basic hypergeometric series. In fact, the monic one-variable little $q$-Jacobi polynomial is independent of $t$ and is given explicitly as
\begin{equation}\label{Lone}
\begin{split}
P_m^L(z)&=\frac{\bigl(qb;q\bigr)_m}{\bigl(q^{m+1}ab;q\bigr)_m(qb)^m}{}_3\phi_2\left(\begin{matrix}
q^{-m},q^{m+1}ab,qbz\\ qb,0\end{matrix};q,q\right)\\
&=\frac{\bigl(qa;q\bigr)_m}{\bigl(q^{m+1}ab;q\bigr)_m}(-1)^mq^{\frac{1}{2}m(m-1)}{}_2\phi_1\left(\begin{matrix}
q^{-m}, q^{m+1}ab\\qa\end{matrix};q,qz\right)
\end{split}
\end{equation}
for $m\in\mathbb{Z}_+$, see  e.g. \cite{KS}. Here the second equality follows from \cite[(III.7)]{GR}.
{}From the second equality in \eqref{Lone} we obtain
\begin{equation}\label{L11}
P_m^L(0)=\frac{\bigl(qa,qab;q\bigr)_m}{\bigl(qab;q\bigr)_{2m}}(-1)^mq^{\frac{1}{2}m(m-1)}.
\end{equation}
{}From the first equality in \eqref{Lone} and the $q$-Vandermonde sum
\begin{equation}\label{qVandermonde}
{}_2\phi_1\left(\begin{matrix} q^{-m},a\\c\end{matrix};q,q\right)=\frac{\bigl(c/a;q\bigr)_m}{\bigl(c;q\bigr)_m}\,a^m
\end{equation}
for $m\in\mathbb{Z}_+$ (see \cite[(II.6)]{GR}), we obtain
\begin{equation}\label{L21}
P_m^L(1)=\frac{\bigl(qb,qab;q\bigr)_m}{\bigl(qab;q\bigr)_{2m}}\bigl(q^ma)^m.
\end{equation}
Finally, from the first equality in \eqref{Lone} we obtain
\begin{equation}\label{L31}
P_m^L(q^{-1}b^{-1})=\frac{\bigl(qb,qab;q\bigr)_m}{\bigl(qab;q\bigr)_{2m}(qb)^m}.
\end{equation}
In the following theorem we give the multivariable analogues of the evaluation formulas
\eqref{L11}, \eqref{L21} and \eqref{L31}. Define $\Delta_\lambda=\Delta_\lambda(a,b;q,t)$
by
\begin{equation}\label{Delta}
\Delta_\lambda(a,b;q,t)=\prod_{1\leq j<k\leq n}
\frac{\bigl(qabt^{2n-j-k+1};q\bigr)_{\lambda_j+\lambda_k}
\bigl(t^{k-j+1};q\bigr)_{\lambda_j-\lambda_k}}
{\bigl(qabt^{2n-j-k};q\bigr)_{\lambda_j+\lambda_k}
\bigl(t^{k-j};q\bigr)_{\lambda_j-\lambda_k}},
\qquad \forall\,\lambda\in\Lambda_n.
\end{equation}
\begin{thm}\label{thmL}
For $\lambda\in\Lambda_n$ we have the evaluation formulas
\begin{equation}\label{L1}
P_\lambda^L(0)=\Delta_\lambda\prod_{i=1}^n\frac{\bigl(qat^{n-i},qabt^{n-i};q\bigr)_{\lambda_i}}{\bigl(qabt^{2(n-i)};q\bigr)_{2\lambda_i}}
(-1)^{\lambda_i}q^{\frac{1}{2}\lambda_i(\lambda_i-1)},
\end{equation}
\begin{equation}\label{L2}
P_\lambda^L(t^\rho)=\Delta_\lambda\prod_{i=1}^n\frac{\bigl(qbt^{n-i}, qabt^{n-i};q\bigr)_{\lambda_i}}
{\bigl(qabt^{2(n-i)};q\bigr)_{2\lambda_i}}\bigl(aq^{\lambda_i}t^{n-i}\bigr)^{\lambda_i},
\end{equation}
\begin{equation}\label{L3}
P_\lambda^L(q^{-1}b^{-1}t^{-\rho})=\Delta_\lambda\prod_{i=1}^n\frac{\bigl(qbt^{n-i},qabt^{n-i};q\bigr)_{\lambda_i}}
{\bigl(qabt^{2(n-i)};q\bigr)_{2\lambda_i}}\bigl(qbt^{n-i}\bigr)^{\lambda_i}.
\end{equation}
\end{thm}
\begin{proof}
The evaluation formulas depend rationally on the parameters, hence it
suffices to prove the theorem for a suitable choice of Zariski dense
complex values of the parameters $a,b,q$ and $t$. The proof now uses
the explicit limit transitions from Macdonald-Koornwinder polynomials to
multivariable little $q$-Jacobi polynomials from \cite{SK} and
\cite{S}. The key ingredient is the following result from the proof
of \cite[Theorem 6.4]{S}. Let $\lambda,\mu\in\Lambda_n$ such that
$\mu<\lambda$. Then
\begin{equation}\label{coefflim}
\lim_{\epsilon\rightarrow 0}\bigl(q^{-\frac{1}{2}}\epsilon\bigr)^{|\lambda|-|\mu|}
c_{\lambda,\mu}(\epsilon^{-1}q^{\frac{1}{2}},-aq^{\frac{1}{2}},\epsilon bq^{\frac{1}{2}}, -q^{\frac{1}{2}};q,t)=
c_{\lambda,\mu}^L(a,b;q,t)
\end{equation}
for the coefficients $c_{\lambda,\mu}$ \eqref{expansionK} and $c_{\lambda,\mu}^L$ \eqref{expansionL}
in the monomial expansions of $P_\lambda$ and $P_\lambda^L$, where $q^{\frac{1}{2}}$ is an arbitrary
choice of square root of $q$.

{\it Proof of \eqref{L1}.}
Observe that
\begin{equation}\label{monlimit1}
\lim_{\epsilon\rightarrow 0}\epsilon^{|\lambda|}m_\lambda(z)=\widetilde{m}_\lambda(0),
\qquad \forall\,\lambda\in\Lambda_n,
\end{equation}
since both sides are equal to zero if $\lambda\in\Lambda_n\setminus\{0\}$
and are equal to one
if $\lambda=0$. Combining \eqref{expansionK}, \eqref{expansionL}, \eqref{coefflim} and
\eqref{monlimit1} now yields
\begin{equation}\label{limpart1}
\lim_{\epsilon\rightarrow 0}\bigl(q^{-\frac{1}{2}}\epsilon\bigr)^{|\lambda|}
P_\lambda\bigl(z;\epsilon^{-1}q^{\frac{1}{2}},-aq^{\frac{1}{2}},\epsilon bq^{\frac{1}{2}}, -q^{\frac{1}{2}};q,t\bigr)=
P_\lambda^L(0;a,b;q,t)
\end{equation}
for $\lambda\in\Lambda_n$. By the evaluation formula \eqref{evalK} for the Macdonald-Koornwinder polynomial $P_\lambda$
and using that the Macdonald-Koornwinder polynomial is symmetric in $a,b,c,d$, we have
\begin{equation*}
\begin{split}
\bigl(q^{-\frac{1}{2}}\epsilon\bigr)^{|\lambda|}
&P_\lambda\bigl(-aq^{\frac{1}{2}}t^\rho;\epsilon^{-1}q^{\frac{1}{2}},
-aq^{\frac{1}{2}},\epsilon bq^{\frac{1}{2}}, -q^{\frac{1}{2}};q,t\bigr)\\
&=\Delta_\lambda(a,b;q,t)
\prod_{i=1}^n\frac{\bigl(qat^{n-i},-qabt^{n-i}\epsilon, qabt^{n-i}, -qat^{n-i}\epsilon^{-1};q\bigr)_{\lambda_i}\epsilon^{\lambda_i}}
{\bigl(qabt^{2(n-i)};q\bigr)_{2\lambda_i}\bigl(-aqt^{n-i}\bigr)^{\lambda_i}}.
\end{split}
\end{equation*}
Combined with the limit transition \eqref{limpart1}, specialized to $z=-aq^{\frac{1}{2}}t^\rho$,
we obtain the evaluation formula \eqref{L1}
for the multivariable little $q$-Jacobi polynomial $P_\lambda^L$.

{\it Proof of \eqref{L2}.}
Now we use the degeneration
\begin{equation}\label{monlimit2}
\lim_{\epsilon\rightarrow 0}\epsilon^{|\lambda|}
m_\lambda(\epsilon^{-1}z)=\widetilde{m}_\lambda(z),\qquad\forall\,\lambda\in\Lambda_n
\end{equation}
for the monomial bases. Together with \eqref{coefflim} this yields the limit transition
\begin{equation}\label{limpart2}
\lim_{\epsilon\rightarrow 0}\bigl(q^{-\frac{1}{2}}\epsilon\bigr)^{|\lambda|}
P_\lambda\bigl(q^{\frac{1}{2}}\epsilon^{-1}z;\epsilon^{-1}q^{\frac{1}{2}},
-aq^{\frac{1}{2}},\epsilon bq^{\frac{1}{2}}, -q^{\frac{1}{2}};q,t\bigr)=
P_\lambda^L(z;a,b;q,t)
\end{equation}
for $\lambda\in\Lambda_n$, which is in accordance with \cite[Theorem 6.4]{S}.
By the evaluation formula \eqref{evalK} for the Macdonald-Koornwinder polynomial $P_\lambda$
we now have
\begin{equation*}
\begin{split}
\bigl(q^{-\frac{1}{2}}\epsilon\bigr)^{|\lambda|}
&P_\lambda\bigl(q^{\frac{1}{2}}\epsilon^{-1}t^\rho;
\epsilon^{-1}q^{\frac{1}{2}},-aq^{\frac{1}{2}},\epsilon bq^{\frac{1}{2}}, -q^{\frac{1}{2}};q,t\bigr)\\
&=\Delta_\lambda(a,b;q,t)\prod_{i=1}^n\frac{\bigl(-qat^{n-i}\epsilon^{-1},-qt^{n-i}\epsilon^{-1}, qbt^{n-i}, qabt^{n-i};q\bigr)_{\lambda_i}\epsilon^{2\lambda_i}}
{\bigl(qabt^{2(n-i)};q\bigr)_{2\lambda_i}\bigl(qt^{n-i}\bigr)^{\lambda_i}}.
\end{split}
\end{equation*}
Combined with the limit transition \eqref{limpart2}, specialized to $z=t^\rho$,
we obtain the evaluation formula \eqref{L2}.

{\it Proof of \eqref{L3}.}
We now use the limit transition
\begin{equation}\label{monlimit3}
\lim_{\epsilon\rightarrow 0}\epsilon^{|\lambda|}m_\lambda(\epsilon z)=\widetilde{m}_\lambda(z^{-1}),\qquad\forall\,
\lambda\in\Lambda_n.
\end{equation}
Together with \eqref{coefflim} this yields
\begin{equation}\label{limpart3}
\lim_{\epsilon\rightarrow 0}\bigl(q^{-\frac{1}{2}}\epsilon\bigr)^{|\lambda|}
P_\lambda\bigl(q^{-\frac{1}{2}}\epsilon z;\epsilon^{-1}q^{\frac{1}{2}},
-aq^{\frac{1}{2}},\epsilon bq^{\frac{1}{2}}, -q^{\frac{1}{2}};q,t\bigr)=
P_\lambda^L(z^{-1};a,b;q,t)
\end{equation}
for $\lambda\in\Lambda_n$.
By the evaluation formula \eqref{evalK} for the Macdonald-Koornwinder polynomial $P_\lambda$
we now have
\begin{equation*}
\begin{split}
\bigl(q^{-\frac{1}{2}}\epsilon\bigr)^{|\lambda|}
&P_\lambda\bigl(\epsilon bq^{\frac{1}{2}}t^\rho;
\epsilon^{-1}q^{\frac{1}{2}},-aq^{\frac{1}{2}},\epsilon bq^{\frac{1}{2}}, -q^{\frac{1}{2}};q,t\bigr)\\
&=\Delta_\lambda(a,b;q,t)\prod_{i=1}^n\frac{\bigl(qbt^{n-i},qabt^{n-i},-qbt^{n-i}\epsilon,-qabt^{n-i}\epsilon;q\bigr)_{\lambda_i}}
{\bigl(qabt^{2(n-i)};q\bigr)_{2\lambda_i}}\bigl(qbt^{n-i}\bigr)^{-\lambda_i}.
\end{split}
\end{equation*}
Combined with the limit transition \eqref{limpart3}, specialized to $z=qbt^\rho$,
we obtain \eqref{L3}.
\end{proof}
The quadratic norms $N_L(\lambda)=N_L(\lambda;a,b;q,t)\in \Fz$ of the multivariable little $q$-Jacobi polynomials
have been explicitly evaluated in \cite{S}. It reads
\[N_L(\lambda)=q^{(\lambda,\lambda)}a^{|\lambda|}t^{2(\rho,\lambda)}
N_L^+(\lambda)N_L^-(\lambda),
\]
where the factors $N_L^{\pm}(\lambda)=N_L^{\pm}(\lambda;a,b;q,t)\in
\Fz$ are given by
\begin{equation*}
\begin{split}
N_L^+(\lambda)&=\prod_{i=1}^n\frac{\bigl(qat^{n-i},qabt^{n-i};q\bigr)_{\lambda_i}}
{\bigl(qabt^{2(n-i)};q\bigr)_{2\lambda_i}}\prod_{1\leq j<k\leq
n}\frac{\bigl(qabt^{2n-j-k+1};q\bigr)_{\lambda_j+\lambda_k}
\bigl(t^{k-j+1};q\bigr)_{\lambda_j-\lambda_k}}
{\bigl(qabt^{2n-j-k};q\bigr)_{\lambda_j+\lambda_k}\bigl(
t^{k-j};q\bigr)_{\lambda_j-\lambda_k}},\\
N_L^-(\lambda)&=\prod_{i=1}^n\frac{\bigl(qt^{n-i},qbt^{n-i};q\bigr)_{\lambda_i}}
{\bigl(q^2abt^{2(n-i)};q\bigr)_{2\lambda_i}}
\prod_{1\leq j<k\leq n}\frac{\bigl(q^2abt^{2n-j-k-1};q\bigr)_{\lambda_j+\lambda_k}
\bigl(qt^{k-j-1};q\bigr)_{\lambda_j-\lambda_k}}
{\bigl(q^2abt^{2n-j-k};q\bigr)_{\lambda_j+\lambda_k}\bigl(
qt^{k-j};q\bigr)_{\lambda_j-\lambda_k}}.
\end{split}
\end{equation*}
We recall from \cite{S} that
\begin{equation}\label{limitLnorms}
\lim_{\epsilon\rightarrow
0}\bigl(q^{-\frac{1}{2}}\epsilon\bigr)^{2|\lambda|}N_K(\lambda;
\epsilon^{-1}q^{\frac{1}{2}},-aq^{\frac{1}{2}},\epsilon bq^{\frac{1}{2}},
-q^{\frac{1}{2}};q,t)=N_L(\lambda;a,b;q,t),\qquad
\forall\,\lambda\in\Lambda_n
\end{equation}
for generic specialized parameters.


\subsection{Multivariable big $q$-Jacobi polynomials}

We repeat the techniques and arguments of the previous subsection to
derive evaluation formulas for multivariable big $q$-Jacobi
polynomials. The monic multivariable big $q$-Jacobi polynomial
\[P_\lambda^B(z)=P_\lambda^B(z;a,b,c,d;q,t)\in\ \Fz[z]^S
\]
of degree
$\lambda\in\Lambda_n$ is of the form
\begin{equation}\label{expansionB}
P_\lambda^B(z)=\widetilde{m}_\lambda(z)+\sum_{\mu\in\Lambda: \mu<\lambda}c_{\lambda,\mu}^B\widetilde{m}_{\mu}(z)
\end{equation}
for certain coefficients $c_{\lambda,\mu}^B=c_{\lambda,\mu}^B(a,b,c,d;q,t)\in \Fz$.
Similarly as for the Koornwinder polynomials
and the multivariable little $q$-Jacobi polynomials, the multivariable big $q$-Jacobi polynomials can be
characterized by a second-order difference equation they should satisfy, or in terms of suitable orthogonality properties,
see \cite{S0}. The characterization as solution of a difference equation is $D_BP_\lambda^B=
E_\lambda^{qJ}P_\lambda^B$ with
\begin{equation*}
\begin{split}
D_B&=\sum_{j=1}^n\bigl(\phi_{B,j}^+(z)(T_j-\hbox{Id})+\phi_{B,j}^-(z)(T_j^{-1}-\hbox{Id})\bigr),\\
\phi_{B,j}^+(z)&=qt^{n-1}\bigl(a-\frac{c}{qz_j}\bigr)\bigl(b+\frac{d}{qz_j}\bigr)\prod_{l\not=j}\frac{z_l-tz_j}{z_l-z_j},\\
\phi_{B,j}^-(z)&=\bigl(1-\frac{c}{z_j}\bigr)\bigl(1+\frac{d}{z_j}\bigr)\prod_{l\not=j}\frac{z_j-tz_l}{z_j-z_l}.
\end{split}
\end{equation*}

We again consider the evaluation formulas for the one-variable big $q$-Jacobi polynomials ($n=1$) first.
The one-variable monic big $q$-Jacobi polynomial is independent of $t$ and is given explicitly as
\begin{equation}\label{Bone}
P_m^B(z)=\frac{\bigl(qa,-qad/c;q\bigr)_m}{\bigl(q^{m+1}ab;q\bigr)_m\bigl(qa/c\bigr)^m}
{}_3\phi_2\left(\begin{matrix} q^{-m},q^{m+1}ab,qaz/c\\ qa,-qad/c \end{matrix};q,q\right)
\end{equation}
for $m\in\mathbb{Z}_+$, see  e.g. \cite{KS}. We immediately obtain
from \eqref{Bone} the evaluation formula
\begin{equation}\label{B11}
P_m^B(c/qa)=\frac{\bigl(qa,-qad/c;q\bigr)_m}{\bigl(q^{m+1}ab;q\bigr)_m}\bigl(c/qa\bigr)^m
\end{equation}
for $m\in\mathbb{Z}_+$.
A straightforward computation using \eqref{Bone} and the $q$-Vandermonde formula
\eqref{qVandermonde} yields the two evaluation formulas
\begin{equation}\label{B21}
\begin{split}
P_m^B(c)&=\frac{\bigl(qa,qab,-qbc/d;q\bigr)_m}{\bigl(qab;q\bigr)_{2m}}d^mq^{\frac{1}{2}m(m-1)},\\
P_m^B(-d)&=\frac{\bigl(-qad/c,qb,qab;q\bigr)_m}{\bigl(qab;q\bigr)_{2m}}(-1)^mc^mq^{\frac{1}{2}m(m-1)}
\end{split}
\end{equation}
for $m\in\mathbb{Z}_+$. Finally, by the $q$-Saalsch{\"u}tz sum
\begin{equation}\label{qSaalschutz}
{}_3\phi_2\left(\begin{matrix}
q^{-m},a,b\\ abq^{1-m}/c,c\end{matrix};
q,q\right)=\frac{\bigl(c/a,c/b;q\bigr)_m}{\bigl(c,c/ab;q\bigr)_m},\qquad
m\in\mathbb{Z}_+
\end{equation}
(see \cite[(II.12)]{GR})
we obtain the evaluation formula
\begin{equation}\label{B31}
P_m^B(-d/qb)=\frac{\bigl(qb,-qbc/d,qab;q\bigr)_m}
{\bigl(qab;q\bigr)_{2m}}\bigl(-d/qb\bigr)^m
\end{equation}
for $m\in\mathbb{Z}_+$.
In the following theorem we give the multivariable analogues of the evaluation formulas
\eqref{B11}, \eqref{B21} and \eqref{B31}. Recall the explicit
expression $\Delta_\lambda$, defined by \eqref{Delta}.
\begin{thm}\label{thmB}
For $\lambda\in\Lambda_n$ we have the evaluation formulas
\begin{equation}\label{B1}
P_\lambda^B(ct^\rho)=\Delta_\lambda\prod_{i=1}^n\frac{\bigl(qat^{n-i},qabt^{n-i},
-\frac{qbc}{d}t^{n-i};q\bigr)_{\lambda_i}}{\bigl(qabt^{2(n-i)};q\bigr)_{2\lambda_i}}
d^{\lambda_i}q^{\frac{1}{2}\lambda_i(\lambda_i-1)},
\end{equation}
\begin{equation}\label{B2}
P_\lambda^B(-dt^\rho)=\Delta_\lambda\prod_{i=1}^n\frac{\bigl(qbt^{n-i},qabt^{n-i},
-\frac{qad}{c}t^{n-i};q\bigr)_{\lambda_i}}
{\bigl(qabt^{2(n-i)};q\bigr)_{2\lambda_i}}(-c)^{\lambda_i}q^{\frac{1}{2}\lambda_i(\lambda_i-1)},
\end{equation}
\begin{equation}\label{B3}
P_\lambda^B\bigl(\frac{c}{qa}\,t^{-\rho}\bigr)=\Delta_\lambda
\prod_{i=1}^n\frac{\bigl(qat^{n-i},qabt^{n-i},-\frac{qad}{c}t^{n-i};q\bigr)_{\lambda_i}}
{\bigl(qabt^{2(n-i)};q\bigr)_{2\lambda_i}}\bigl(\frac{c}{qa}t^{i-n}\bigr)^{\lambda_i},
\end{equation}
\begin{equation}\label{B4}
P_\lambda^B\bigl(-\frac{d}{qb}\,t^{-\rho}\bigr)=\Delta_\lambda
\prod_{i=1}^n\frac{\bigl(qbt^{n-i},qabt^{n-i},-\frac{qbc}{d}t^{n-i};q\bigr)_{\lambda_i}}
{\bigl(qabt^{2(n-i)};q\bigr)_{2\lambda_i}}\bigl(-\frac{d}{qb}t^{i-n}\bigr)^{\lambda_i}.
\end{equation}
\end{thm}
\begin{proof}
The proof is analogous to the proof of Theorem \ref{thmB}, so we only indicate the main steps.
{}From (the proof of) \cite[Theorem 7.5]{S} (replacing the role of the limiting parameter $\epsilon$
in \cite[Theorem 7.5]{S} by $\bigl(cd/q\bigr)^{\frac{1}{2}}\epsilon$), we obtain the two limit transitions
\begin{equation}\label{limitB1}
\lim_{\epsilon\rightarrow 0}\epsilon^{|\lambda|}P_\lambda\Bigl(\epsilon^{-1}z;
\frac{qa\epsilon}{c}, -\frac{qb\epsilon}{d}, c\epsilon^{-1},
-d\epsilon^{-1};q,t\Bigr)=P_\lambda^B(z;a,b,c,d;q,t),
\end{equation}
\begin{equation}\label{limitB2}
\lim_{\epsilon\rightarrow 0}\epsilon^{|\lambda|}P_\lambda\Bigl(\epsilon z;
\frac{qa\epsilon}{c}, -\frac{qb\epsilon}{d}, c\epsilon^{-1},
-d\epsilon^{-1};q,t\Bigr)=P_\lambda^B(z^{-1};a,b,c,d;q,t)
\end{equation}
for $\lambda\in\Lambda_n$ and for a specific Zariski-dense subset of
specialized parameters $a,b,c,d,q$ and $t$. If we set
$z=\epsilon^{-1}ct^\rho$ (respectively $z=-\epsilon^{-1}dt^\rho$) in
\eqref{limitB1} and use the evaluation formula \eqref{evalK} for the
Macdonald-Koornwinder polynomial $P_\lambda$, then we arrive at the
evaluation formula \eqref{B1} (respectively \eqref{B2}). The
evaluation formulas \eqref{B3} and \eqref{B4} follow from
\eqref{evalK} and the limit transition \eqref{limitB2} specialized
to $z=\frac{qa}{c}\,t^\rho$ and $z=-\frac{qb}{d}\,t^\rho$,
respectively.
\end{proof}
We end this section by recalling the quadratic norm evaluations of
the multivariable big $q$-Jacobi polynomials from \cite{S}.
We define an $\Fz$-linear functional
\[h_B=h_B^{a,b,c,d;q,t}: \Fz[z]^S\rightarrow \Fz
\]
by requiring
\begin{equation*}
h_B\bigl(P_\lambda^B\bigr)=
\begin{cases}
1,\qquad &\hbox{ if } \lambda=0,\\
0,\qquad &\hbox{ if } \lambda\not=0.
\end{cases}
\end{equation*}
An analytic definition of $h_B$ can be given for specific Zariski
dense choice of specialized parameters as a multidimensional
$q$-integral with explicit weight function, see e.g. \cite{S}. We
define a sesqui-linear form $\langle \cdot,\cdot\rangle_B= \langle
\cdot,\cdot\rangle_B^{a,b,c,d;q,t}$ on $\Fz[z]^S$ by
\begin{equation}\label{Bpairing}
\langle
p_1,p_2\rangle_B=h_B\bigl(p_1(z)\overline{p_2(z)}\bigr),\qquad
p_1,p_2\in \Fz[z]^S,
\end{equation}
with the same convention on the complex conjugation as for the
multivariable little $q$-Jacobi polynomials. The orthogonality relations for the
multivariable big $q$-Jacobi polynomials are now given by
\begin{equation}
\langle P_\lambda^B,P_\mu^B\rangle_B=N_B(\lambda)\delta_{\lambda,\mu},\qquad
\forall\,\lambda,\mu\in\Lambda_n
\end{equation}
with the quadratic norms $N_B(\lambda)=N_B(\lambda;a,b,c,d;q,t)\in \Fz$
given by
\[N_B(\lambda)=(cd)^{|\lambda|}t^{(\rho,\lambda)}
\left(\prod_{i=1}^nq^{{\lambda_i\choose
2}}\Bigl(-\frac{qbc}{d}t^{n-i},-\frac{qad}{c}t^{n-i};q\Bigr)_{\lambda_i}\right)N_L^+(\lambda)N_L^-(\lambda),
\]
with $N_L^{\pm}=N_L^{\pm}(a,b;q,t)$ the same factors as for the
multivariable little $q$-Jacobi polynomials. We recall from \cite{S}
that
\begin{equation}\label{limitBnorms}
\lim_{\epsilon\rightarrow 0}\epsilon^{2|\lambda|}N_K
\bigl(\lambda;\frac{qa\epsilon}{c}, -\frac{qb\epsilon}{d}, c\epsilon^{-1},
-d\epsilon^{-1};q,t\bigr)=N_B(\lambda;a,b,c,d;q,t),\qquad
\lambda\in\Lambda_n,
\end{equation}
for generic specialized parameters.


\section{Generalized dimension formulas}\label{Dimensionsection}

In this section we consider the generalized dimensions
$D_q(\lambda;a,b;t)$ (see \eqref{Dgeneralized}) in detail. We
suppress the dependence on the parameters $a,b$ and $t$ as much as
possible, so we write $D_q(\lambda)=D_q(\lambda;a,b;t)$. We will
show in Section \ref{quantum} that $D_q(\lambda)$ for specialized
parameters gives the quantum dimension for the irreducible spherical
representation $V_\lambda^q$ ($\lambda\in\Lambda_n$) associated to
the standard quantum complex Grassmannian. The degeneration ($q=1$)
of the quantum dimensions yields the complex dimensions of the
irreducible $\Le$-spherical $\Ko$-representations associated to the
complex and real Grassmannian, see Subsection \ref{cG}. The
degeneration ($q=0$), which corresponds to the dimensions of the
irreducible $\LF$-spherical $\KF$-representations $V_\lambda^{\Fb}$
for $\Fb$ non-Archimedean (see Subsection \ref{repp}) is analyzed in
detail in Subsection \ref{Genpdim}.

\subsection{Generalized quantum dimensions}

We begin by giving an explicit expression for the generalized dimension
$D_q(\lambda)$.

\begin{lem}\label{explicitD}
We have
\[D_q(\lambda)=a^{-|\lambda|}t^{-2(\rho,\lambda)}\prod_{i=1}^nv_i(\lambda_i;q)
\prod_{1\leq j<k\leq n}w_{j,k}^+(\lambda_j+\lambda_k;q)w_{j,k}^-(\lambda_j-\lambda_k;q),\qquad \lambda\in\Lambda_n,
\]
where $v_i(m;q)$, $w_{j,k}^{\pm}(m;q)$ for $m\in\mathbb{Z}_+$ are given by
\begin{equation*}
\begin{split}
v_i(m;q)&=\frac{\bigl(at^{n-i}, q^{-1}abt^{n-i};q\bigr)_m}
{\bigl(qt^{n-i},bt^{n-i};q\bigr)_m}\frac{(1-q^{2m-1}abt^{2(n-i)})}
{(1-q^{-1}abt^{2(n-i)})},\\
w_{j,k}^+(m;q)&=\frac{\bigl(q^{-1}abt^{2n-j-k+1};q\bigr)_m}
{\bigl(abt^{2n-j-k-1};q\bigr)_m}\frac{(1-q^{m-1}abt^{2n-j-k})}
{(1-q^{-1}abt^{2n-j-k})},\\
w_{j,k}^-(m;q)&=\frac{\bigl(t^{k-j+1};q\bigr)_m}{\bigl(qt^{k-j-1};q\bigr)_m}
\frac{(1-q^mt^{k-j})}{(1-t^{k-j})}.
\end{split}
\end{equation*}
\end{lem}
\begin{proof}
This follows from a direct computation using the evaluation formulas and the quadratic
norm formulas for the multivariable little $q$-Jacobi polynomials.
\end{proof}
The explicit expressions for the generalized dimensions associated
to the fundamental partitions $\omega_r\in\Lambda_n$
($r=0,\ldots,n$) simplify drastically.
\begin{prop}\label{lemref}
For $r=0,\ldots,n$ we have
\begin{equation}\label{fundamentalqformula}
D_q(\omega_r)=\frac{\bigl(qabt^{2n-r-1}, t^{n+1-r}, at^{n-r}, abt^{2n-r};t\bigr)_r}
{\bigl(q, t, bt^{n-r}, abt^{n-r-1};t\bigr)_r}
\left(\frac{1-abt^{2n-2r-1}}{1-abt^{2n-1}}\right)a^{-r}t^{r(r+1-2n)}.
\end{equation}
\end{prop}
\begin{proof}
Note that
\[|\omega_r|=r,\qquad 2(\rho,\omega_r)=r(2n-r-1).
\]
By Lemma \ref{explicitD} we thus obtain
\begin{equation}\label{Dq}
D_q(\omega_r)=a^{-r}t^{r(r+1-2n)}\prod_{i=1}^rv_i(1;q)
\prod_{1\leq j<k\leq r}w_{j,k}^+(2;q)\prod_{j=1}^r\prod_{k=r+1}^nw_{j,k}^+(1;q)
\prod_{j=1}^r\prod_{k=r+1}^nw_{j,k}^-(1;q).
\end{equation}
By direct computations we obtain the explicit expressions
\begin{equation}\label{form}
\begin{split}
\prod_{i=1}^rv_i(1;q)&=\frac{\bigl(at^{n-r}, q^{-1}abt^{n-r};t\bigr)_r\bigl(qabt^{2n-2r};t^2\bigr)_r}
{\bigl(qt^{n-r}, bt^{n-r};t\bigr)_r\bigl(q^{-1}abt^{2n-2r};t^2\bigr)_r},\\
\prod_{1\leq j<k\leq r}w_{j,k}^+(2;q)&=\frac{\bigl(abt^{2n-r};t\bigr)_r\bigl(q^{-1}abt^{2n-2r}, qabt^{2n-2r-1};t^2\bigr)_r}
{\bigl(q^{-1}abt^{2n-2r}, qabt^{2n-2r-1}, abt^{2n-2r-1};t\bigr)_r}
\frac{(1-abt^{2n-2r-1})}{(1-abt^{2n-1})},\\
\prod_{j=1}^r\prod_{k=r+1}^nw_{j,k}^+(1;q)&=\frac{\bigl(q^{-1}abt^{2n-2r}, abt^{2n-2r-1};t\bigr)_r}
{\bigl(q^{-1}abt^{n-r}, abt^{n-r-1};t\bigr)_r},\\
\prod_{j=1}^r\prod_{k=r+1}^nw_{j,k}^-(1;q)&=\frac{\bigl(t^{n+1-r},
qt^{n-r};t\bigr)_r}{\bigl(t,q;t\bigr)_r}.
\end{split}
\end{equation}
Only the second formula needs explanation. Set
\[e_r(u;t)=\prod_{1\leq j<k\leq
r}\left(\frac{1-ut^{1-j-k}}{1-ut^{-j-k}}\right),
\]
then we can write
\[\prod_{1\leq j<k\leq r}w_{j,k}^+(2;q)=e_r(q^{-1}abt^{2n})
e_r(abt^{2n})e_r(abt^{2n-1})e_r(qabt^{2n-1}).
\]
On the other hand, by induction to $r\in\{1,\ldots,n\}$,
\[e_r(u;t)=\frac{\bigl(ut^{2-2r};t^2\bigr)_{r-1}}{\bigl(ut^{1-2r};t\bigr)_{r-1}}.
\]
Combining both formulas we easily obtain the desired expression for
$\prod_{1\leq j<k\leq r}w_{j,k}^+(2;q)$. The lemma now follows by
substituting the expressions \eqref{form} in \eqref{Dq}, and by
simplifying the resulting expression.
\end{proof}

As a curiosity, note that the elementary identity
\[\bigl(ut^{1-r};t\bigr)_r=(-u)^rt^{-\frac{1}{2}r(r-1)}\bigl(u^{-1};t\bigr)_r
\]
allows us to rewrite \eqref{fundamentalqformula} as
\[
D_q\bigl(\omega_r)=\frac{\bigl(t^{-n},a^{-1}t^{1-n},
a^{-1}b^{-1}t^{1-2n},
q^{-1}a^{-1}b^{-1}t^{2-2n};t\bigr)_r}
{\bigl(t,q,b^{-1}t^{1-n},
a^{-1}b^{-1}t^{2-n};t\bigr)_r}\\
\left(\frac{1-a^{-1}b^{-1}t^{1+2r-2n}}{1-a^{-1}b^{-1}t^{1-2n}}\right)
\bigl(qa t^{2n-1}\bigr)^r
\]
for $r=0,\ldots,n$, which one recognizes as the weight function for $t$-Racah polynomials,
see e.g. \cite{KS}.


\subsection{Generalized $p$-adic dimensions}\label{Genpdim}

In this subsection we consider the $p$-adic degeneration of the
generalized quantum dimension $D_q(\lambda)$ and we relate it to the
complex dimension of the $\LF$-spherical irreducible
$\KF$-representation $V_\lambda^{\mathbb{F}}$ for non-Archimedean
local fields $\mathbb{F}$ (see Section \ref{local}). Denote
$\delta_{k,l}$ ($k,l\in\mathbb{Z}$) for the Kronecker delta
function: it is one if $k=l$ and zero otherwise.
\begin{lem}\label{specialq=0}
The functions $v_i(m;q)$, $w_{j,k}^{\pm}(m;q)$ \textup{(}see Lemma \ref{explicitD}\textup{)}
for $m\in\mathbb{Z}_+$ are regular at $q=0$.
If we write $v_i(m)=v_i(m;0)$ and $w_{j,k}^{\pm}(m)=w_{j,k}^{\pm}(m;0)$ for their constant terms,
then $v_i(0)=w_{j,k}^{\pm}(0)=1$ and
\begin{equation*}
\begin{split}
v_i(m)&=\left(\frac{1-at^{n-i}}{1-bt^{n-i}}\right)t^{i-n}\bigl(1-(1-\delta_{m,1})abt^{n-i}\bigr),\\
w_{j,k}^{+}(m)&=\left(\frac{1-abt^{2n-j-k+1-\delta_{m,1}}}{1-abt^{2n-j-k-1}}\right)t,\\
w_{j,k}^{-}(m)&=\frac{1-t^{k-j+1}}{1-t^{k-j}}
\end{split}
\end{equation*}
for $m\geq 1$.
\end{lem}
\begin{proof}
As an example, we compute $v_i(m)$ for $m\in\mathbb{Z}_{\geq 1}$ (the other computations
are similar). To compute $v_i(1)$ we rewrite $v_i(1;q)$ as
\[v_i(1;q)=\left(\frac{1-qabt^{2(n-i)}}{1-qt^{n-i}}\right)\left(\frac{1-at^{n-i}}{1-bt^{n-i}}\right)
\left(\frac{q-abt^{n-i}}{q-abt^{2(n-i)}}\right),
\]
which yields
\[v_i(1)=v_i(1;0)=\left(\frac{1-at^{n-i}}{1-bt^{n-i}}\right)t^{i-n}.
\]
To compute $v_i(m)$ for $m\geq 2$ we write $v_i(m;q)$ as
\begin{equation*}
\begin{split}
v_i(m;q)=&\left(\frac{(1-q^{2m-1}abt^{2(n-i)})\bigl(qat^{n-i};q\bigr)_{m-1}\bigl(qabt^{n-i};q\bigr)_{m-2}}
{\bigl(qt^{n-i};q\bigr)_m\bigl(qbt^{n-i};q\bigr)_{m-1}}\right)\\
&\times \left(\frac{(1-at^{n-i})(1-abt^{n-i})}{(1-bt^{n-i})}\right)
\left(\frac{q-abt^{n-i}}{q-abt^{2(n-i)}}\right),
\end{split}
\end{equation*}
which yields
\[v_i(m)=v_i(m;0)=\left(\frac{1-at^{n-i}}{1-bt^{n-i}}\right)t^{i-n}(1-abt^{n-i}),\qquad m\geq 2.
\]
\end{proof}
\begin{prop}\label{formula}
Let $\lambda\in\Lambda_n$. The generalized quantum dimension $D_q(\lambda)$ \textup{(}see \eqref{Dgeneralized}\textup{)}
is regular at $q=0$.
Writing $D_0(\lambda)=D_0(\lambda;a,b;t)$ for $D_q(\lambda;a,b;t)|_{q=0}$, we have
\[D_0(\lambda)=a^{-|\lambda|}t^{-2(\rho,\lambda)}\prod_{i=1}^nv_i(\lambda_i)\prod_{1\leq j<k\leq n}
w_{j,k}^+(\lambda_j+\lambda_k)w_{j,k}^-(\lambda_j-\lambda_k),\qquad \lambda\in\Lambda_n.
\]
\end{prop}
\begin{proof}
This is immediate from Lemma \ref{explicitD} and the previous lemma.
\end{proof}
\begin{eg}
For $n=1$, so that $\Lambda_1=\mathbb{Z}_+$, we have $D_0(0)=1$ and
\[D_0(\lambda)=a^{-\lambda}\left(\frac{1-a}{1-b}\right)\bigl(1-(1-\delta_{\lambda,1})ab\bigr),\qquad \lambda\geq 1.
\]
\end{eg}
In contrast to the quantum case, the explicit expression for
$D_0(\lambda)=D_0(\lambda;a,b;t)$ from Proposition \ref{formula}
simplifies for all partitions $\lambda\in\Lambda_n$ as follows.

\begin{thm}\label{padicformula}
For $\lambda\in\Lambda_n$ we have
\begin{equation}\label{padicfor}
D_0(\lambda)=a^{-|\lambda|}t^{-2(\rho,\lambda)}\bigg[{n \atop \partial\lambda'}
\bigg]_t \frac{\bigl(at^{n-\lambda_1'};t\bigr)_{\lambda_1'}\bigl(abt^{2n-\lambda_1'-\lambda_2'};t\bigr)_{\lambda_1'+\lambda_2'}}
{\bigl(bt^{n-\lambda_1'}, abt^{n-\lambda_1'-1};t\bigr)_{\lambda_1'}}
\left(\frac{1-abt^{2n-2\lambda_1'-1}}{1-abt^{2n-1}}\right).
\end{equation}
\end{thm}

As remarked already in Subsection \ref{repp}, the generalized
$p$-adic dimension formula \eqref{padicfor} combined with its
representation theoretic interpretation \eqref{reldimD0} from
\cite{OU1} leads to the explicit dimension formulas
\eqref{padic-dimension} for the irreducible $\LF$-spherical
$\KF$-representations $V_\lambda^{\Fb}$ ($\lambda\in\Lambda_n$).

\begin{cor}\label{fundamentalcase}
For $r=0,\ldots,n$ we have
\begin{equation}\label{Dpadic}
D_0(\omega_r)=
\frac{\bigl(t^{n+1-r}, at^{n-r}, abt^{2n-r};t\bigr)_r}
{\bigl(t, bt^{n-r}, abt^{n-r-1};t\bigr)_r}
\left(\frac{1-abt^{2n-2r-1}}{1-abt^{2n-1}}\right)a^{-r}t^{r(r+1-2n)}.
\end{equation}
\end{cor}
\begin{proof}
Follows either by specializing Theorem \ref{padicformula} to
$\lambda=\omega_r$, or by taking $q=0$ in
\eqref{fundamentalqformula}.
\end{proof}

We next proceed to prove Theorem \ref{padicformula}, which is based
on the expression for $D_0(\lambda)$ as given in Proposition \ref{formula}.
We divide the proof into several elementary lemmas. Since the theorem is obviously
correct for $\lambda=0$, we fix for the proof a nonzero partition
$\lambda\in\Lambda_n$.

\begin{lem}\label{el1}
\[
\prod_{i=1}^n v_{i}(\lambda_i)=t^{\frac{1}{2}\lambda_1'(\lambda_1'+1)-n\lambda_1'}
\frac{\bigl(at^{n-\lambda_1'};t\bigr)_{\lambda_1'}
\bigl(abt^{n-\lambda_2'};t\bigr)_{\lambda_2'}}{\bigl(bt^{n-\lambda_1'};t\bigr)_{\lambda_1'}}.
\]
\end{lem}
\begin{proof}
Using the explicit expression for $v_i$ (see Lemma \ref{specialq=0}), we
have
\[
\begin{split}
\prod_{i=1}^n v_{i}(\lambda_i)&=\prod_{i=1}^{\lambda_2'}
\left(\frac{1-at^{n-i}}{1-bt^{n-i}}\right)t^{i-n}\bigl(1-abt^{n-i}\bigr)
\prod_{i=\lambda_2'+1}^{\lambda_1'}\left(\frac{1-at^{n-i}}{1-bt^{n-i}}\right)t^{i-n}\\
&=\prod_{i=1}^{\lambda_1'} \left(\frac{1-at^{n-i}}{1-bt^{n-i}}\right)t^{i-n}
\prod_{i=1}^{\lambda_2'}\left({1-abt^{n-i}}\right)\\
&=t^{\frac{1}{2}\lambda_1'(\lambda_1'+1)-n\lambda_1'}\frac{\bigl(at^{n-\lambda_1'};t\bigr)_{\lambda_1'}
\bigl(abt^{n-\lambda_2'};t\bigr)_{\lambda_2'}}{\bigl(bt^{n-\lambda_1'};t\bigr)_{\lambda_1'}},
\end{split}
\]
as desired.
\end{proof}

\begin{lem}\label{el2}
\[
\prod_{1\leq j<k\leq n} w_{j,k}^-(\lambda_j-\lambda_k)=\bigg[{n
\atop \partial\lambda'}\bigg]_t.
\]
\end{lem}
\begin{proof} Since $w_{j,k}^-(0)=1$ and $w_{j,k}^-(m)=\frac{1-t^{k-j+1}}{1-t^{k-j}}$
for $m\ne 0$ (see Lemma \ref{specialq=0}), the term $w_{j,k}^{-}(\lambda_j-\lambda_k)$
contributes to the product only if $\lambda_j \ne \lambda_k$.
Hence
\[\prod_{1\leq j<k\leq n} w_{j,k}^-(\lambda_j-\lambda_k)=A/B,
\]
with
\[
\begin{split}
&A=\prod_{1\leq j<k\leq
n}\frac{1-t^{k-j+1}}{1-t^{k-j}}=\prod_{l=1}^{n}\frac{(1-t^{l+1})^{n-l}}{(1-t^{l})^{n-l}}=[n]_t! \\
&B=\prod_{\substack{1\leq j<k\leq n \\
\lambda_j =
\lambda_k}}\frac{1-t^{k-j+1}}{1-t^{k-j}}=\prod_{l\ge
0}\prod_{1\leq j<k\leq \lambda_l'-\lambda_{l+1}'}\frac{1-t^{k-j+1}}{1-t^{k-j}}=\prod_{l \ge 0}
[\lambda_l'-\lambda_{l+1}']_t!,
\end{split}
\]
as desired.
\end{proof}

\begin{lem}\label{el3}
\[
\prod_{1\leq j<k\leq n} w_{j,k}^+(\lambda_j+\lambda_k)=t^{\frac{1}{2}\lambda_1'(\lambda_1'-1)+\lambda_1'(n-\lambda_1')}
\frac{(abt^{2n-\lambda_1'};t)_{\lambda_1'-1}}{(abt^{2n-2\lambda_1'};t)_{\lambda_1'-1}}
\frac{(abt^{2n-2\lambda_1'-1};t)_{\lambda_1'}(abt^{2n-\lambda_1'-\lambda_2'};t)_{\lambda_2'}}
{(abt^{n-\lambda_1'-1};t)_{\lambda_1'}(abt^{n-\lambda_2'};t)_{\lambda_2'}}.
\]
\end{lem}
\begin{proof} The term $w_{j,k}^+(m)$ has three possible values depending on $m$ being $0$, $1$ or $\ge 2$
(see Lemma \ref{specialq=0}).
We are hence led to consider the following sets
\[
\begin{split}
 X_1&=\{(j,k) \,|\, 1\leq j<k\leq n,\,\, \lambda_j+\lambda_k=1\}=
 \{(j,k) \,|\, \lambda_2' < j \le \lambda_1' < k \le n \}\\
X_2&=\{(j,k) \,|\, 1\leq j<k\leq n,\,\,  \lambda_j+\lambda_k \ge 2 \}=X_2^a\sqcup X_2^b,
\end{split}
\]
where the disjoint subsets $X_2^a$ and $X_2^b$ are given by
\[X_2^a=\{(j,k) \,|\, 1\leq j<k\leq \lambda_1'\},\qquad
X_2^b=\{(j,k) \,|\, 1\leq j\leq \lambda_2',\,\, \lambda_1'<k\leq n\}.
\]
We then have
\[\prod_{1\leq j<k\leq n} w_{j,k}^+(\lambda_j+\lambda_k) =
\Pi(X_1)\Pi(X_2^a)\Pi(X_2^b)
\]
where we write $\Pi(X)=\prod_{(j,k)\in X} w_{j,k}^+(\lambda_j+\lambda_k)$
for a subset $X\subseteq \{(j,k) \, | \, 1\leq j<k\leq n \}$.
We have
\[
\begin{split}
\Pi(X_1)&=t^{(n-\lambda_1')(\lambda_1'-\lambda_2')}\prod_{j=\lambda_2'+1}^{\lambda_1'}\prod_{k=\lambda_1'+1}^n
\left(\frac{1-abt^{2n-j-k}}{1-abt^{2n-j-k-1}}\right)\\
&=t^{(n-\lambda_1')(\lambda_1'-\lambda_2')}\prod_{j=\lambda_2'+1}^{\lambda_1'}
\left(\frac{1-abt^{2n-j-\lambda_1'-1}}{1-abt^{n-j-1}}\right)= t^{(n-\lambda_1')(\lambda_1'-\lambda_2')}
\frac{\bigl(abt^{2n-2\lambda_1'-1};t\bigr)_{\lambda_1'-\lambda_2'}}{\bigl(ab
t^{n-\lambda_1'-1};t\bigr)_{\lambda_1'-\lambda_2'}}.
\end{split}
\]
Using the same notations and simplification arguments as in the
proof of Proposition \ref{lemref}, we have
\[\Pi(X_2^a)=t^{\frac{1}{2}\lambda_1'(\lambda_1'-1)}e_{\lambda_1'}(abt^{2n})
e_{\lambda_1'}(abt^{2n-1})=
t^{\frac{1}{2}\lambda_1'(\lambda_1'-1)}\frac{\bigl(abt^{2n-\lambda_1'};t\bigr)_{\lambda_1'-1}}
{\bigl(abt^{2n-2\lambda_1'};t\bigr)_{\lambda_1'-1}}.
\]
By a similar computation as for $\Pi(X_1)$, we finally have
\[\Pi(X_2^b)=t^{(n-\lambda_1')\lambda_2'}\,
\frac{\bigl(abt^{2n-\lambda_1'-\lambda_2'},
abt^{2n-\lambda_1'-\lambda_2'-1};t\bigr)_{\lambda_2'}}
{\bigl(abt^{n-\lambda_2'}, abt^{n-\lambda_2'-1};t\bigr)_{\lambda_2'}}.
\]
Combining these expressions leads to the desired result.
\end{proof}
We can now use Proposition \ref{formula}, Lemma \ref{el1}, Lemma \ref{el2} and Lemma
\ref{el3} to obtain an explicit expression for the generalized
$p$-adic dimension $D_0(\lambda)$. After straightforward
simplifications, we arrive at the desired expression \eqref{padicfor} for
$D_0(\lambda)$.


\section{The one-parameter family of quantum complex Grassmannians}\label{quantum}

In this section we consider the one-parameter family of
quantum complex Grassmannians from \cite{NDS}, \cite{DS} and
\cite{OS} (see also \cite{L} for the general theory of harmonic
analysis on quantum compact symmetric spaces). This more general
set-up does no longer fit into the framework of quantum Gelfand pairs
associated to \textsf{CQG} algebras as discussed in Section \ref{CQG},
since the role of the quantum subgroup is replaced by a suitable infinitesimal
counterpart, depending on an additional continuous parameter.

In \cite{NDS} and \cite{DS} harmonic analysis on the one-parameter
family of quantum complex Grassmannians was studied in detail,
leading to the interpretation of Macdonald-Koornwinder polynomials
and multivariable big and little $q$-Jacobi polynomials as the associated quantum
zonal spherical functions. {}From these results we now derive
representation theoretic interpretations of the
evaluation formulas and norm formulas for Macdonald-Koornwinder polynomials and multivariable big and
little $q$-Jacobi polynomials. The arguments follow closely \cite[Section 6]{N},
now applied to the one-parameter family of quantum complex Grassmannians.
In this more general context we discuss the quantum dimension
formulas for irreducible spherical representations in terms of
quantum zonal spherical functions. The analysis leads to the
quantum dimension formulas \eqref{qcf} and \eqref{qcfcons} for the
irreducible $\mathcal{U}_q$-representations $V_\lambda^q$ ($\lambda\in\Lambda_n$).

We now first recall the construction of the quantum zonal spherical functions on a
one-parameter family of quantum analogues of the complex
Grassmannian $\Ko/\Le$. We will freely use the notations from Section
\ref{qGrassmannian}.
The role of the subgroup $\Le$ is
taken over by a two-sided co-ideal $\mathfrak{k}_\sigma\subseteq \mathcal{U}_q$
depending (in a suitable sense continuously) on an additional parameter $\sigma\in\mathbb{R}\cup \{\infty\}$, see e.g.
\cite{NDS}, \cite{DS} for the explicit definition of $\mathfrak{k}_\sigma$.
The case $\sigma=\infty$ corresponds to the standard quantum
complex Grassmannian from Section \ref{qGrassmannian}.

The space of $\mathfrak{k}_\sigma$-invariants in $\RR_q(\Ko)$,
\[\RR_q^\sigma(\Ko/\Le):=\{f\in\RR_q(\Ko) \,\, | \,\, f(YX)=0,\qquad \forall\, X\in\mathfrak{k}_\sigma,\,\,
\forall\, Y\in\mathcal{U}_q \}
\]
is a $\mathcal{U}_q$-submodule of $\RR_q(\Ko)$ with respect to the regular $\mathcal{U}_q$-action
\[\bigl(X\,f\bigr)(X^\prime)=f(S(X)X^\prime),\qquad f\in\RR_q^\sigma(\Ko/\Le),\quad X,X^\prime\in\mathcal{U}_q.
\]
The $\mathcal{U}_q$-module $\RR_q^\sigma(\Ko/\Le)$ serves as a
one-parameter family of quantum analogues of the Grassmann
representation $\RR(\Ko/\Le)$. Furthermore, $\RR_q^\infty(\Ko/\Le)$
is the Grassmann representation $\RR_q(\Ko/\Le)$ associated to the
standard quantum complex Grassmannian from Section
\ref{qGrassmannian}. The space
\[\bigl(L_\mu^q\bigr)^{\mathfrak{k}_\sigma}=\{v\in L_\mu^q \,\, | \,\, X v=0,\quad \forall\, X\in\mathfrak{k}_\sigma \}
\]
for $\mu\in P_d^+$ serves as the one-parameter quantum analogue of the space of $\Le$-fixed vectors
in $L_\mu^{\mathbb{C}}$. The subspace $\bigl(L_\mu^q\bigr)^{\mathfrak{k}_\sigma}$ is
one-dimensional or zero dimensional, and it is one-dimensional if and
only if $\mu\in\Lambda_n^\natural$. Consequently, the irreducible
spherical representations $V_\lambda^q=L_{\lambda^\natural}^q$
($\lambda\in\Lambda_n$) associated to the standard quantum complex
Grassmannian also serve as complete set of representatives of
finite dimensional, irreducible $\mathcal{U}_q$-representations
with nonzero (in fact, one-dimensional) space of
$\mathfrak{k}_\sigma$-invariants.
{}From the Peter-Weyl decomposition for $\RR_q(\Ko)$ we obtain the irreducible decomposition
\[\RR_q^\sigma(\Ko/\Le)\simeq \bigoplus_{\lambda\in\Lambda_n}
V_{\lambda}^q
\]
as $\mathcal{U}_q$-modules. Furthermore, for $\mu\in P_d^+$ we
recall from \cite{DS} that
\[\bigl(L_\mu^q\bigr)^{\mathfrak{k}_\infty}=\bigl(L_\mu^q\bigr)^{\mathcal{U}_q(\mathfrak{m}_{\mathbb{C}})},
\]
cf. Subsection \ref{qGr}.

We fix generic $\sigma,\tau\in\mathbb{R}$ unless specified differently (the case $\sigma=\infty$ and/or $\tau=\infty$
will subsequently be analyzed as degeneration of the theory for $\sigma, \tau \in\mathbb{R}$).
For $\lambda\in\Lambda_n$ we choose a nonzero $\mathfrak{k}_\sigma$-fixed vector
$v_\sigma(\lambda)\in \bigl(V_\lambda^q\bigr)^{\mathfrak{k}_\sigma}$ which we normalize
by requiring
\[v_\sigma(\lambda)=\sum_{\nu\in\mathbb{Z}^{d}}v_\sigma(\lambda;\nu),\qquad
v_\sigma(\lambda;\lambda^\natural)=u_{\lambda^\natural},
\]
where $v_\sigma(\lambda;\nu)\in V_\lambda^q$
is the weight $\nu$ component of $v_\sigma(\lambda)$ in
$V_\lambda^q$ and $u_{\lambda^\natural}$ is the highest weight vector of
$V_\lambda^q=L_{\lambda^\natural}^q$ (it is known that $v_\sigma(\lambda)$ has a weight $\lambda^\natural$ component,
hence such a normalization is possible). {}From the analysis in \cite[Section 6]{DS} it furthermore
follows that $\bigl(V_\lambda^q\bigr)^{\mathfrak{k}_\sigma^*}$ is spanned by $K^{-\delta}v_{\sigma+d-2n}(\lambda)$
for $\lambda\in\Lambda_n$ (see also \cite{OS} for the special case $d=2n$).

The normalized quantum zonal spherical functions
$\varphi_\lambda^{\sigma,\tau}$ ($\sigma,\tau\in\mathbb{R}$) associated to the one-parameter family
of quantum Grassmannians are now defined by
\begin{equation}\label{zsf}
\varphi_\lambda^{\sigma,\tau}(\,\cdot\,)=q^{\langle \delta,\lambda^\natural\rangle}
\frac{\langle \cdot\, v_\sigma(\lambda),
K^{-\delta}v_{\tau+d-2n}(\lambda)\rangle_{\lambda}}
{\langle u_{\lambda^\natural},u_{\lambda^\natural}\rangle_{\lambda}},\qquad \lambda\in\Lambda_n.
\end{equation}
The $\varphi_\lambda^{\sigma,\tau}$
($\lambda\in\Lambda_n$) form a basis of
\[\mathcal{H}_{\sigma,\tau}:=\{f\in \RR_q(\Ko) \,\, | \,\, f(YZ)=0=f(ZX),\quad
Z\in\mathcal{U}_q,\,\, X\in \mathfrak{k}_\sigma,\,\,
Y\in\mathfrak{k}_\tau \}.
\]
The analogue of the quantum dimension formula \eqref{dimforquantum} in terms of the quantum zonal
spherical functions $\varphi_\lambda^{\sigma,\tau}$ ($\lambda\in\Lambda_n$) reads as follows.

\begin{lem}
For $\lambda\in\Lambda_n$ and generic $\sigma,\tau\in\mathbb{R}$ we have
\begin{equation}\label{dimforquantumst}
\textup{Dim}_q\bigl(V_\lambda^q\bigr)=\frac{\varphi_\lambda^{\sigma,\sigma+2n-d}(K^{\delta})\varphi_\lambda^{\tau+d-2n,\tau}(K^\delta)}
{\|\varphi_\lambda^{\sigma,\tau}\|_h^2}.
\end{equation}
\end{lem}
\begin{proof}
By the quantum Schur orthogonality relations (see Proposition \ref{Schur}{\bf b}) we have
\begin{equation}\label{ato}
\|\varphi_\lambda^{\sigma,\tau}\|_h^2=\frac{q^{2\langle\delta,\lambda^\natural\rangle}}{\textup{Dim}_q(V_\lambda^q)}
\frac{\|v_\sigma(\lambda)\|_\lambda^2 \|v_{\tau+d-2n}(\lambda)\|_\lambda^2}{\|u_{\lambda^\natural}\|_\lambda^4},
\end{equation}
where $\|v\|_{\lambda}^2:=\langle v, v\rangle_{\lambda}$. On the other hand, from the definition \eqref{zsf} of
the quantum zonal spherical function $\varphi_\lambda^{\sigma,\tau}$
we have
\begin{equation}\label{bto}
\varphi_\lambda^{\sigma,\sigma+2n-d}(K^\delta)=\frac{\|v_\sigma(\lambda)\|_\lambda^2}{\|u_{\lambda^\natural}\|_\lambda^2}\,
q^{\langle\delta,\lambda^\natural\rangle}.
\end{equation}
The desired quantum dimension formula follows directly from \eqref{ato} and \eqref{bto}.
\end{proof}

To arrive at an explicit evaluation of the quantum dimensions, we now translate \eqref{dimforquantumst} into an expression involving the evaluations
and the quadratic norms of Macdonald-Koornwinder polynomials.

The radial part $P_\lambda^{\sigma,\tau}:=\hbox{Res}_T\bigl(\varphi_\lambda^{\sigma,\tau}\bigr)\in
\mathbb{C}[u^{\pm 1}]$ is defined to be the unique Laurent polynomial in $d$
independent variables $u_1,\ldots,u_{d}$ such that
\[P_\lambda^{\sigma,\tau}(q^\nu)=
\varphi_\lambda^{\sigma,\tau}(K^{\nu})\qquad\forall\,\nu\in\mathbb{Z}^{d}.
\]
We have the following key facts from \cite{NDS}, see also \cite{DS} and \cite{OS}.
Set $z_j=u_ju_{d+1-j}^{-1}$ ($j=1,\ldots,n$), then
\[P_\lambda^{\sigma,\tau}\in \mathbb{C}[z^{\pm 1}]^W,\qquad
\forall\lambda\in \Lambda_n,
\]
where $W=S_n\ltimes\{\pm 1\}^n$ as before. {}From now on we will
view $P_\lambda^{\sigma,\tau}$ as
$W$-invariant Laurent polynomial in the independent variables
$z_1,\ldots,z_n$. The key result is the identification of $P_\lambda^{\sigma,\tau}$
with Macdonald-Koornwinder polynomials,
\[P_\lambda^{\sigma,\tau}(z)=
P_\lambda(z;-q^{\sigma+\tau+1},-q^{-\sigma-\tau+1},q^{\sigma-\tau+1},
q^{-\sigma+\tau+2(d-2n)+1};q^2,q^2),\qquad\lambda\in\Lambda_n
\]
(the particular normalization of $\varphi_\lambda^{\sigma,\tau}$
ensures that it precisely coincides with the monic
Macdonald-Koornwinder polynomial).
We denote $\langle \cdot,\cdot\rangle_{\sigma,\tau}$ for the corresponding orthogonality
pairing \eqref{Kpairing}, with the six parameters $(a,b,c,d,q,t)$
specialized to
\[(-q^{\sigma+\tau+1},-q^{-\sigma-\tau+1},q^{\sigma-\tau+1},
q^{-\sigma+\tau+2(d-2n)+1},q^2,q^2).
\]
We set $\|p\|_{\sigma,\tau}^2=\langle p,p\rangle_{\sigma,\tau}^2$
for the corresponding quadratic norm of $p\in\mathbb{C}[z^{\pm 1}]^W$.

\begin{lem}
For $\lambda\in\Lambda_n$ and generic $\sigma,\tau\in\mathbb{R}$ we have
\begin{equation*}
\begin{split}
\varphi_\lambda^{\sigma,\tau}(K^\delta)&=P_\lambda^{\sigma,\tau}(q^{d-2n+1}q^{2\rho}),\\
\|\varphi_{\lambda}^{\sigma,\tau}\|_h^2&=
\|P_\lambda^{\sigma,\tau}\|_{\sigma,\tau}^2.
\end{split}
\end{equation*}
In particular,
\begin{equation}\label{normevaluationcombi}
\hbox{Dim}_q\bigl(V_\lambda^q\bigr)=\frac{P_\lambda^{\sigma,\sigma+2n-d}(q^{d-2n+1}q^{2\rho})
P_\lambda^{\tau+d-2n,\tau}(q^{d-2n+1}q^{2\rho})}
{\|P_\lambda^{\sigma,\tau}\|_{\sigma,\tau}^2}.
\end{equation}
\end{lem}
\begin{proof}
Observe that $u=q^\delta\in \bigl(\mathbb{R}^\times\bigr)^d$ maps to
$z=q^{d-2n+1}q^{2\rho}\in\bigl(\mathbb{R}^\times\bigr)^n$ under the assignment
$z_j=u_ju_{d+1-j}^{-1}$ ($j=1,\ldots,n$). By the definition of $P_\lambda^{\sigma,\tau}$ as the radial
part of $\varphi_\lambda^{\sigma,\tau}$, we obtain the desired expression of $\varphi_\lambda^{\sigma,\tau}(K^\delta)$
in terms of Macdonald-Koornwinder polynomials. The quadratic norm formula follows from \cite[Corollary 8.1]{DS}.
The quantum dimension formula \eqref{dimforquantumst} then yields \eqref{normevaluationcombi}.
\end{proof}
We are now in a position to derive the explicit quantum dimension formulas for the spherical
$\mathcal{U}_q$-representations $V_\lambda^q$ ($\lambda\in\Lambda_n$).
\begin{cor}\label{quantumcomplexformulacor}
Proposition \ref{quantumcomplexformula} holds.
\end{cor}
\begin{proof}
Observe that the value of the Macdonald-Koornwinder polynomial
\[P_\lambda^{\sigma,\sigma+2n-d}(z)=
P_{\lambda}(z;-q^{2\sigma+2n-d+1}, -q^{-2\sigma-2n+d+1},
q^{d-2n+1}, q^{d-2n+1};q^2,q^2),
\]
at $z=q^{d-2n+1}q^{2\rho}$  can be evaluated in closed form
by the evaluation formula \eqref{evalK} for Macdonald-Koornwinder polynomials.
The resulting expression is
\begin{equation*}
\begin{split}
P_\lambda^{\sigma, \sigma+2n-d}(q^{d-2n+1}q^{2\rho})
&=q^{(2n-d-1)|\lambda|-2\langle\rho,\lambda\rangle}\\
&\times\prod_{1\leq
j<k\leq
n}\frac{(1-q^{2(d-2n+1+\lambda_j+\lambda_k+\rho_j+\rho_k)})
(1-q^{2(\lambda_j-\lambda_k
+\rho_j-\rho_k)})}
{(1-q^{2(d-2n+1+\rho_j+\rho_k)})(1-q^{2(\rho_j-\rho_k)})}\\
\times&\prod_{i=1}^n\frac{\bigl(q^{2(d-2n+1+\rho_i)},q^{2(d-2n+1+\rho_i)},
-q^{2(\sigma+1+\rho_i)},
-q^{2(-\sigma+d-2n+1+\rho_i)};q^2\bigr)_{\lambda_i}}
{\bigl(q^{2(d-2n+1+2\rho_i)};q^2\bigr)_{2\lambda_i}}.
\end{split}
\end{equation*}
On the other hand, by the closed expression \eqref{normKoorn}
for the quadratic norms of the Macdonald-Koornwinder polynomials we have
\begin{equation*}
\begin{split}
\|P_\lambda^{\sigma,\tau}\|_{\sigma,\tau}^2&=\prod_{j=1}^n
\left\{\frac{\bigl(q^{2(1+\rho_j)},
q^{2(d-2n+1+\rho_j)};q^2\bigr)_{\lambda_j}^2}{\bigl(q^{2(d-2n+1+2\rho_j)},
q^{2(d-2n+2+2\rho_j)};q^2\bigr)_{2\lambda_j}}\right.\\
&\left.\times\bigl(-q^{2(\sigma+1+\rho_j)},
-q^{2(-\sigma+d-2n+1+\rho_j)}, -q^{2(\tau+d-2n+1+\rho_j)},
-q^{2(-\tau+1+\rho_j)};q^2\bigr)_{\lambda_j}\right\}.
\end{split}
\end{equation*}
Combining these explicit formulas with \eqref{normevaluationcombi} we obtain an explicit expression
for $\hbox{Dim}_q(V_\lambda^q)$, which results in Proposition \ref{quantumcomplexformula}.
\end{proof}

We end this section by considering the $\sigma\rightarrow \infty$
and $\sigma=\tau\rightarrow\infty$ degenerations in the above
representation theoretic formulas. This in particular leads to the
interpretation of $D_q(\lambda;a,b;t)$ for special values of the
parameters as quantum dimensions of $V_\lambda^q$
($\lambda\in\Lambda_n$), see \eqref{qcfcons}.

In \cite{DS} the limit $\sigma\rightarrow\infty$ (respectively $\sigma,\tau\rightarrow\infty$)
is considered, leading to the interpretation of the
multivariable big $q$-Jacobi polynomials $P_\lambda^B(\cdot;1,q^{2(d-2n)},1,q^{2\tau+2(d-2n)};q^2,q^2)$
(respectively the multivariable little $q$-Jacobi polynomials $P_\lambda^L(\cdot;q^{2(d-2n)},1;q^2,q^2)$)
as radial parts of quantum zonal spherical functions associated to the standard quantum Grassmann
representation $\RR_q^\infty(\Ko/\Le)=\RR_q(\Ko/\Le)$.
In these degenerations, the quantum dimensions $\hbox{Dim}_q(V_\lambda^q)$
are naturally expressed in terms of the multivariable big and little $q$-Jacobi polynomials as follows.

\begin{prop}\label{propdim}
We have
\begin{equation*}
\textup{Dim}_q(V_\lambda^q)=
\bigl(-q^{2\tau+d-2n-1}\bigr)^{|\lambda|}\frac{P_\lambda^L(0;q^{2(d-2n)},1;q^2,q^2)P_\lambda^{\tau+d-2n,\tau}(q^{d-2n+1}q^{2\rho})}
{N_B(\lambda;1,q^{2(d-2n)},1,q^{2\tau+2(d-2n)};q^2,q^2)},\qquad \forall\,\lambda\in\Lambda_n.
\end{equation*}
Furthermore, in terms of the generalized quantum dimension $D_q(\lambda;a,b;t)$,
see \eqref{Dgeneralized}, we have
\[\textup{Dim}_q(V_\lambda^q)=
\frac{\bigl(P_\lambda^L(0;q^{2(d-2n)},1;q^2,q^2)\bigr)^2}{N_L(\lambda;q^{2(d-2n)},1;q^2,q^2)}=
D_{q^2}(\lambda;q^{2(d-2n+1)}, q^2;q^2),\qquad \forall\,\lambda\in\Lambda_n.
\]
\end{prop}
\begin{proof}
We degenerate the expression \eqref{normevaluationcombi} of the quantum dimension $\hbox{Dim}_q(V_\lambda^q)$
by taking the limit $\sigma\rightarrow\infty$ (respectively $\sigma=\tau\rightarrow\infty$).
In order to parallel the particular degeneration scheme from \cite{DS}, we make use of the elementary symmetry
\begin{equation}\label{symmKo}
P_\lambda(-z;-a,-b,-c,-d;q,t)=(-1)^{|\lambda|}P_\lambda(z;a,b,c,d;q,t)
\end{equation}
for the monic Macdonald-Koornwinder polynomials, which implies the symmetry
\begin{equation}\label{symmnorm}
N_K(\lambda;-a,-b,-c,-d;q,t)=N_K(\lambda;a,b,c,d;q,t),\qquad\lambda\in\Lambda_n
\end{equation}
for the corresponding quadratic norms (which also immediately follows from its explicit evaluation \eqref{evalK}).
Combined with the limit formula \eqref{limitBnorms} in base $q^2$, with parameters $(a,b,c,d,t)$ specialized to
$(1,q^{2(d-2n)},1,q^{2\tau+2(d-2n)},q^2)$ and with $\epsilon=q^{\sigma+\tau-1}$, we obtain
\[\lim_{\sigma\rightarrow\infty}q^{2(\sigma+\tau-1)|\lambda|}\|P_\lambda^{\sigma,\tau}\|_{\sigma,\tau}^2=
N_B(\lambda;1,q^{2(d-2n)},1,q^{2\tau+2(d-2n)};q^2,q^2),\qquad \forall\,\lambda\in\Lambda_n.
\]
On the other hand, from \eqref{L1} and the explicit evaluation formula from the proof of Corollary \ref{quantumcomplexformulacor}
we obtain
\[\lim_{\sigma\rightarrow\infty}\bigl(-q^{2\sigma}\bigr)^{|\lambda|}P_\lambda^{\sigma,\sigma+2n-d}(q^{d-2n+1}q^{2\rho})=
q^{(d-2n+1)|\lambda|}P_\lambda^L(0;q^{2(d-2n)},1;q^2;q^2),\qquad \forall\,\lambda\in\Lambda_n.
\]
Taking the limit $\sigma\rightarrow\infty$ in \eqref{normevaluationcombi} now leads to the expression of the quantum dimension
as mixture of factors involving Macdonald-Koornwinder, multivariable big and multivariable little $q$-Jacobi polynomials.

Taking $\tau=\sigma$ in \eqref{normevaluationcombi} and using the fact that
\[\lim_{\sigma\rightarrow\infty}\bigl(q^{-1+2\sigma}\bigr)^{2|\lambda|}\|P_\lambda^{\sigma,\sigma}\|_{\sigma,\sigma}^2=
N_L(\lambda;q^{2(d-2n)},1;q^2,q^2),\qquad
\forall\,\lambda\in\Lambda_n
\]
in view of \eqref{limitLnorms}, we similarly obtain the desired expression of the quantum dimension
$\hbox{Dim}_q(V_\lambda^q)$ in terms of $D_q(\lambda;a,b;t)$.
\end{proof}

\begin{rema}\label{standardrema}
By the arguments from \cite{DS} one can show that
\[\lim_{\sigma\rightarrow\infty}\left(\frac{\varphi_\lambda^{\sigma,\sigma+2n-d}}{\|\varphi_{\lambda}^{\sigma,\sigma+2n-d}\|_h}\right)=
\varphi_\lambda^q,\qquad \lambda\in\Lambda_n
\]
in the pre-Hilbert space $\RR_q(\Ko)$, where $\varphi_\lambda^q$ ($\lambda\in\Lambda_n$) are the (orthonormal) quantum zonal spherical
functions associated to the standard quantum complex Grassmanian. Degenerating \eqref{dimforquantumst}
accordingly, we obtain the quantum dimension formula \eqref{dimforquantum}.
\end{rema}


\end{document}